\documentclass[a4paper,fleqn]{cas-sc}
\usepackage[numbers]{natbib}
\usepackage{ulem}
\usepackage{placeins}
\usepackage{amsthm}

\theoremstyle{remark}
\newtheorem{remark}{Remark}
\usepackage{bm}
\usepackage{tikz}
\usetikzlibrary{arrows.meta,calc,positioning}
\usepackage{url}
\makeatletter
\providecommand\captionof[1]{\def\@captype{#1}\caption}
\makeatother
\ExplSyntaxOn
\keys_define:nn { cas / fig }
  { unknown .code:n = { \tl_if_blank:nT {#1} { \tl_set:Nx \l_fig_pos_tl { \l_keys_key_str } } } }
\keys_define:nn { cas / tbl }
  { unknown .code:n = { \tl_if_blank:nT {#1} { \tl_set:Nx \l_tbl_pos_tl { \l_keys_key_str } } } }
\tl_put_right:Nn \__reset_fig: { \tl_set:Nn \l_fig_pos_tl { tbp } }
\tl_put_right:Nn \__reset_tbl: { \tl_set:Nn \l_tbl_pos_tl { tbp } }
\ExplSyntaxOff
\DeclareRobustCommand{\_}{\textunderscore\allowbreak}
\ExplSyntaxOn
\RenewDocumentCommand \printorcid { }
  {
    \seq_if_empty:NF \g_stm_orcid_seq
      {
        \group_begin:
          \tex_let:D \thefootnote \relax \footnotetext
            { \raggedright \textsc{orcid}(s):\c_space_token \seq_use:Nn \g_stm_orcid_seq { ;~ } }
        \group_end:
      }
  }
\ExplSyntaxOff

\newcommand{\ii}{\mathrm{i}}

\newcommand{\pml}{\mathrm{pml}}
\newcommand{\phys}{\mathrm{phys}}

\providecommand{\upDelta}{\mathrm{\Delta}}
\newcommand{\laplace}{\mathop{}\!{\upDelta}}

\newcommand{\range}{\operatorname{range}}
\newcommand{\argmin}{\operatorname*{arg\,min}}
\newcommand{\mat}[1]{\boldsymbol{#1}}
\newcommand{\vect}[1]{\boldsymbol{#1}}
\newcommand{\experimentlabel}[1]{%
  \par\addvspace{5pt}\noindent\textbf{#1.}\enspace\ignorespaces}
\newcounter{algorithm}
\newenvironment{paperalgorithm}[2]
{%
  \refstepcounter{algorithm}%
  \par\smallskip\noindent
  \begin{minipage}{0.94\textwidth}
  \small
  \hrule
  \vspace{2pt}
  \noindent\textbf{Algorithm~\thealgorithm. #1}\label{#2}
  \vspace{2pt}
  \hrule
  \vspace{3pt}
}
{%
  \vspace{2pt}
  \hrule
  \end{minipage}
  \par\smallskip
}

\newlength{\panelunit}

\renewcommand{\appendix}{\par
  \setcounter{section}{0}%
  \setcounter{subsection}{0}%
  \renewcommand{\thesection}{\Alph{section}}}

\ExplSyntaxOn
\cs_set:Npn \__first_footerline: {\small\itshape Preprint}
\ExplSyntaxOff
\begin{document}
\let\WriteBookmarks\relax
\def\floatpagepagefraction{1}
\def\textpagefraction{.001}

\shorttitle{Trajectory-aware DD-MSRB preconditioner for parametric Helmholtz problems}
\shortauthors{A. Ping, Z. Peng, P. Yang}

\title[mode=title]{A Trajectory-Aware Domain-Decomposed Multispace Reduced Basis Preconditioner for Parametric Helmholtz Problems}

\author[1]{An Ping}
\author[1]{Zhichao Peng}
\cormark[1]
\ead{pengzhic@ust.hk}
\cortext[1]{Corresponding author}
\author[2]{Pengliang Yang}
\affiliation[1]{organization={Department of Mathematics, The Hong Kong University of Science and Technology},
  city={Hong Kong}, country={China}}
\affiliation[2]{organization={Laoshan Laboratory},
  city={Qingdao}, country={China}}

\begin{abstract}
Repeated solutions of parametric Helmholtz problems arise in multi-query applications such as sensitivity analysis and uncertainty quantification. Classical preconditioned iterative solvers do not automatically reuse low-dimensional structures across parameters. Multispace reduced basis (MSRB) preconditioners can exploit such structures by constructing iteration-specific reduced order models (ROMs) from training data, but storing global iteration-dependent bases incurs substantial memory costs. Building on trajectory-aware MSRB training, we develop a domain-decomposed MSRB (DD-MSRB) preconditioner for parametric Helmholtz systems truncated by perfectly matched layers (PMLs). Within flexible GMRES, each reduced space is trained on inverse-action snapshots along the Krylov trajectories generated by the  ROM-based preconditioners applied in preceding iterations. To reduce storage, problem-adapted spatial partitions exploit heterogeneous local snapshot compressibility across the physical interior and the PML. Local basis dimensions in each subdomain are determined greedily under a global snapshot-error budget to balance approximation accuracy and storage. The resulting hybrid solver applies DD-MSRB preconditioning during the initial iterations to capture dominant cross-parameter modes and then switches, when necessary, to a shifted-Laplacian preconditioner to reduce error components poorly represented by the reduced spaces. Numerical experiments on two-dimensional half-space propagation, elliptic-obstacle scattering, and three-dimensional waveguides demonstrate significant acceleration over the shifted-Laplacian baseline and substantial memory savings relative to global MSRB. For the largest waveguide system, with 11.0 million unknowns, DD-MSRB reduces deployed ROM storage by 44.9\% and mean online solve time by 30.2\% compared with global MSRB.
\end{abstract}

\begin{keywords}
Helmholtz equation \sep domain-decomposed reduced order model \sep preconditioner \sep perfectly matched layer
\MSC[2020]{65N22, 65N35, 65F10}
\end{keywords}
\maketitle
\section{Introduction}
\label{sec:intro}
The Helmholtz equation describes time-harmonic wave propagation in applications
including acoustic and seismic imaging, photonic devices, and waveguide
analysis. Many such problems are posed in unbounded domains, where an outgoing
radiation condition selects the physically relevant solution. Numerical
computation requires a finite domain with a truncation that absorbs outgoing
waves with minimal reflection~\cite{BaylissTurkel1980}. Perfectly matched layers (PMLs) are widely applied for this
purpose~\cite{Berenger1994}. PML truncation and spatial discretization lead to
large indefinite systems, with complex coordinate stretching generally making
the operator non-Hermitian. A large number of degrees of freedom (DOFs) may be
required to resolve the oscillatory solution and control the pollution
error~\cite{BabuskaSauter1997}, making the storage requirements of classical
direct methods substantial in 3D. Moreover, standard Krylov methods may
stagnate without effective
preconditioning~\cite{Erlangga2008,ErnstGander2012}.

A substantial literature has been devoted to making a single Helmholtz solve
more efficient. Established approaches to preconditioning Helmholtz systems include
shifted-Laplacian
preconditioners~\cite{BaylissGoldsteinTurkel1983,ErlanggaVuikOosterlee2004,ErlanggaOosterleeVuik2006,VanGijzenErlanggaVuik2007,Gander2006}, domain-decomposition
methods~\cite{BENAMOU199768,farhat2000two,GanderMagoulesNataf2002,GrahamSpenceVainikko2017}, sweeping
preconditioners~\cite{EngquistYing2011,taus2020sweeps},
deflation~\cite{DwarkaVuik2022}, and time-domain preconditioners~\cite{bristeau1998controllability,stolk2021time,grote2019controllability,appelo2020waveholtz}.
We refer to the overview papers~\cite{LahayeTangVuik2017,GanderZhang2019,GanderZhang2022} for a more comprehensive literature review.

In multi-query applications such as inverse problems, design optimization, and uncertainty quantification, Helmholtz systems are solved repeatedly for various parameter values describing material or geometric configurations. Repeated solves at varying parameters raise the complementary question of how
information generated by previous solves can be exploited to accelerate subsequent
solves. Reduced order models (ROMs)~\cite{BennerGugercinWillcox2015} provide a
strategy to use such information through an offline--online decomposition.
The offline stage constructs the ROM by extracting low-dimensional structures
from available solution data for a parametric problem, while the online stage
exploits the ROM to efficiently predict solutions at new parameter values
through projection or interpolation~\cite{BennerGugercinWillcox2015}.
Such ROMs have been developed for parametric Helmholtz and time-harmonic Maxwell
systems in~\cite{SenVeroyHuynhDeparisNguyenPatera2006,ChenHesthavenMadayRodriguez2010,pomplun2010accelerated,hess2013fast,benner2015uncertainty}. 

The cheap low-fidelity approximations provided by ROMs can be used as initial guesses and preconditioners to accelerate full-order Krylov solvers that compute high-fidelity solutions.
ROM-based preconditioners have been developed for parametric PDEs~\cite{dal2017numerical,DalSanto2018,nguyen2018reduced,liu2022rom,mcclarren2022data,peng2024reduced,Peng2025}. Among these, the multispace reduced basis (MSRB) preconditioner is particularly attractive because it constructs a separate reduced space for each Krylov iteration, accounting for the circular dependence between the iteration trajectory and the ROM-based preconditioner. A recent variant tailored to transport sweeps for the parametric radiative transfer equation (RTE)~\cite{TangPeng2025} showed that standard ROM-based preconditioners lose effectiveness as iterations proceed if they fail to account for how preconditioners in earlier steps change the iteration trajectories used to construct subsequent reduced spaces. In summary, ROMs for subsequent iterations should be constructed in a trajectory-aware manner using data
generated by preceding iterations with their corresponding ROM-based
preconditioners.

Extending the trajectory-aware MSRB framework to large-scale, high-frequency Helmholtz--PML systems incurs substantial storage costs: resolving oscillatory wave fields and padding the domain with a PML require many DOFs, while storing a distinct reduced basis at each iteration further increases the storage requirement. Localized and domain-decomposed model reduction offers a promising remedy by allowing smaller local bases in more compressible regions. Far-field reduction~\cite{BuffoniTelibIollo2009} exploits regional differences
in solution complexity. Domain-decomposed ROMs construct local reduced basis or least-squares
Petrov--Galerkin approximations and couple them across
subdomains~\cite{BaigesCodinaIdelsohn2013,HoangChoiCarlberg2021}.
Component-based formulations include the reduced basis element
method~\cite{MadayRonquist2002} and static-condensation methods for
Helmholtz acoustics~\cite{HuynhKnezevicPatera2013Complex}, with port
reduction further compressing interface
spaces~\cite{EftangPatera2013,smetana2016optimal}.
Local approximation spaces can also be enriched
adaptively~\cite{ohlberger2015error,carlberg2015adaptive}
or constructed using randomized
sampling~\cite{BuhrSmetana2018}.
A comprehensive review of localized ROMs is given in~\cite{BuhrEtAl2021}.

To build an efficient ROM-based preconditioner for parametric Helmholtz equations while mitigating its storage overhead, we propose a hybrid domain-decomposed multispace reduced basis (DD-MSRB) preconditioner. We construct domain-decomposed ROMs for the early Krylov iterations to capture dominant parametric modes with reduced storage, and then switch to a shifted-Laplacian preconditioner to damp the
remaining error modes not well represented by the ROMs.  The main contributions
of this work are as follows.
\begin{enumerate}
\item We develop a trajectory-aware, hybrid MSRB-enhanced preconditioner for
Helmholtz--PML systems within the flexible generalized minimal residual method (FGMRES), extending our previous work on the RTE~\cite{TangPeng2025}. The design exploits the complementary strengths of MSRB and classical shifted-Laplacian preconditioners: ROM-based preconditioners rapidly damp
error components associated with the dominant solution structures in early
iterations, while a  shifted-Laplacian preconditioner more
effectively damps the remaining error modes that are poorly represented
by the ROMs.

\item To reduce basis storage, especially the overhead associated
with PML unknowns, we introduce a problem-adapted DD-MSRB
construction that exploits variations in local snapshot
compressibility across oscillatory interior fields, far-field
components, and rapidly decaying PML profiles. A greedy algorithm
selects the local basis dimensions under a global snapshot-error 
budget to balance compression accuracy and storage.

\item For a three-dimensional parametric waveguide problem, we empirically
observe that the domain-decomposed ROMs better capture dominant solution
structures in subdomain interiors than global ROMs, despite potentially
larger errors near subdomain interfaces. We observe that these interface errors are
effectively damped by the shifted-Laplacian preconditioner. Compared with
global MSRB, the resulting hybrid DD-MSRB preconditioner converges faster
after switching to the shifted-Laplacian preconditioner.
\end{enumerate}

The paper is organized as follows.  Section~\ref{sec:problem} introduces the
parametric Helmholtz system and the shifted-Laplacian preconditioner.  Section~\ref{sec:dd-rom} presents the domain-decomposed proper orthogonal decomposition (POD) for constructing spatially local ROMs, together with a local rank selection to balance memory and accuracy.
Section~\ref{sec:correction-rom} develops trajectory-aware training and the
hybrid preconditioner.  Section~\ref{sec:numerics} demonstrates the performance of the proposed methods through a series of 2D and 3D tests.
Section~\ref{sec:conclusion} concludes the paper.

\section{Background}
\label{sec:problem}

Wave propagation in unbounded media is governed by the exterior Helmholtz
problem.  In multi-query applications such as optimal design and control,
inverse problems, and uncertainty quantification, this equation
must be solved repeatedly for parameter values that represent material
properties, source configurations, or geometric features.  This leads naturally
to parametric Helmholtz problems.

\subsection{Parametric Helmholtz--PML systems}
\label{subsec:helmholtz-pml}
In this paper, we consider the following
parametric Helmholtz equation
\begin{equation}
  \laplace u(x;\mu)
  + \frac{\omega^2}{c(x;\mu)^2}\,u(x;\mu)
  = f(x;\mu),
\label{eq:helmholtz-phys}
\end{equation}
supplemented by the Sommerfeld radiation condition at infinity.  Here, $u$ is the scalar wave field, $\omega>0$ is the
angular frequency, $c(x;\mu)$ is the wavespeed, $f$ is the source, and $\mu$ denotes the parameters that determine the material, source, or geometric configuration.
In practice, the open domain is truncated to a bounded
physical domain of interest $\Omega$ by surrounding it with a perfectly matched layer (PML)
$\Omega_{\pml}$~\cite{Berenger1994}.  The resulting computational domain is
$\Omega_h=\Omega\cup\Omega_{\pml}$.  Let $\sigma_x,\sigma_y\geq0$ be damping
profiles supported in $\Omega_{\pml}$, and define the complex stretching
factors
\begin{equation}
  s_x(x)=1+\frac{\ii\sigma_x(x)}{\omega},
  \qquad
  s_y(y)=1+\frac{\ii\sigma_y(y)}{\omega}.
\label{eq:stretching}
\end{equation}
In all experiments, we use the quadratic profile
\(\sigma_x(x)=\sigma_{\max}(\xi_x(x)/L)^2\), where \(\xi_x(x)\) is the
distance from \(x\) to the inner boundary of the PML and \(L\) is the PML
thickness; \(\sigma_y\) (and \(\sigma_z\) in three dimensions) is defined
analogously.
The corresponding PML operator is
\begin{equation}
  \mathcal L_{\pml}(\mu)u
  =
  \frac{1}{s_x}\partial_x\!\left(\frac{1}{s_x}\partial_x u\right)
  +
  \frac{1}{s_y}\partial_y\!\left(\frac{1}{s_y}\partial_y u\right)
  +
  \frac{\omega^2}{c(x;\mu)^2}\,u ,
\label{eq:pml-operator}
\end{equation}
and the truncated problem reads
\begin{equation}
  \mathcal L_{\pml}(\mu)u=f_e
  \qquad\text{in }\Omega_h ,
\label{eq:pml-continuous}
\end{equation}
where $f_e$ denotes the extension of the physical source to the extended
domain. A homogeneous Dirichlet condition is imposed at the outer PML boundary. After a  second-order central finite-difference
discretization of this operator, we obtain the linear system
\(\mat A_{\mu} \vect u=\vect b_{\mu}\).

\subsection{Shifted-Laplacian framework for Helmholtz systems}
\label{subsec:helmholtz-solvers}
Shifted-Laplacian preconditioning~\cite{ErlanggaOosterleeVuik2006,VanGijzenErlanggaVuik2007} is a well-established and competitive
approach for solving indefinite Helmholtz systems and is therefore adopted
as the baseline iterative solver in this work.
Here, we introduce additional complex damping into an auxiliary
Helmholtz operator.
With the sign convention in \eqref{eq:helmholtz-phys}, let
\[
 \mat M_\beta(\mu)=\mat A_{\mu}+\ii\beta\mat D_{\mu},
 \qquad
 \mat D_{\mu}=\operatorname{diag}\!\left(\frac{\omega^2}{c(x_i;\mu)^2}\right),
 \qquad \beta>0.
\]
The shifted mass term leaves the original Helmholtz system unchanged and
introduces additional absorption in the preconditioning operator, so that its inverse is easier to approximate
by incomplete factorizations or multigrid. 

With \(\mathcal B_\beta\approx\mat M_\beta^{-1}\), right-preconditioned
GMRES solves
\[
 \mat A_{\mu}\mathcal B_\beta\vect y=\vect b,
 \qquad \vect u=\mathcal B_\beta\vect y.
\]
Thus, \(\mathcal B_\beta\) serves as an approximate inverse of \(\mat A_\mu\), where the shift factor $\beta$ governs the trade-off between efficiency and effectiveness: stronger damping makes the auxiliary system easier to solve but makes its approximation error of the original Helmholtz operator larger~\cite{Gander2006,CocquetGander2017}.

The approximate inverse \(\mathcal B_\beta\) of \(\mat M_\beta\) can be implemented by incomplete
LU factorization (ILU) or multigrid (MG). In this work, we use
ILU for the two-dimensional problems, which lead to moderate-size linear systems, because of its better efficiency in this regime, and MG for the larger-scale three-dimensional problems for its better scalability.

\subsection{Motivation: reusable inverse actions and storage-efficient representations}
\label{subsec:motivation-strategy}

Shifted-Laplacian preconditioners are effective for standalone solves, but they treat each parameter in isolation and thus do not exploit shared information across solves.
In multi-query applications, similar systems are repeatedly solved over a nearby parameter range, and their inverse actions and Krylov directions often share common low-dimensional structures.
To exploit such structures in new parameter solves, MSRB preconditioners~\cite{DalSanto2018} build a sequence of ROM-based preconditioners tailored to individual Krylov iterations.

However, these iteration-dependent ROMs with global basis vectors can incur substantial memory overhead in Helmholtz--PML systems.
Inverse-action snapshots typically exhibit strongly heterogeneous compressibility across the computational domain, for example, oscillating in the  interior and rapidly decaying  within the PML.
A global ROM nevertheless stores every retained mode over the entire PML-extended domain, wasting storage in highly compressible regions.
To overcome this bottleneck, we introduce a domain-decomposed representation that constructs local reduced spaces and allocates their ranks adaptively according to regional compressibility.

\section{Domain-decomposed reduced order models\label{sec:dd-rom}}
In preparation for the MSRB preconditioners, we introduce a domain-decomposed POD (DD-POD) to exploit the heterogeneous compressibility of snapshots across the computational domain.  The offline stage builds local bases under a global snapshot-error budget, with local ranks in each subdomain determined by a greedy selection. The online stage applies Galerkin
projection at a new parameter (Section~\ref{sec:dd-rom-online}).
\subsection{Offline basis construction\label{sec:dd-rom-offline}}
\label{subsec:rom-projection}
\label{subsec:domain-decomposed-rom}
For a parametric linear system \(\mat A_h(\mu)\vect u_h(\mu)=\vect f_h(\mu)\),
classical POD collects snapshots \(\mat E=[\vect u_h(\mu_1),\dots,\vect u_h(\mu_{n_s})]\) and
constructs a reduced basis \(\mat \Phi_r\in\mathbb{C}^{N_h\times r}\) via a
truncated singular value decomposition (SVD)
\(\mat E\approx\mat \Phi_r\mat \Sigma_r\mat V_r^*\).

For a $d$-dimensional PML of relative thickness $\theta$ on each side, the fraction of PML DOFs is $1-\frac{1}{(1+2\theta)^d}$. If $\theta=0.15$, the PML fraction becomes $40.8\%$ in two dimensions and $54.5\%$ in three dimensions. Applying a global POD basis in the entire computational domain is inefficient when solutions are highly compressible in subregions such as the PML.

\textbf{Local basis construction.} To address this issue, we replace the global POD basis by local POD bases on
nonoverlapping subdomains \(\Omega_b\), \(b=1,\dots,B\), satisfying
\(\Omega_h=\bigcup_{b=1}^B \Omega_b\).  Let \(n_b\) be the number of DOFs in
\(\Omega_b\) satisfying $\sum_{b=1}^B n_b=N_h$, and let
\(\mathcal{I}_b=(i_{1}^{(b)},\dots,i_{n_b}^{(b)})\) be the corresponding index
set.  The restriction matrix \(\mat R_b\) and its zero-padding adjoint
\(\mat R_b^*\) satisfy
  \begin{equation*}
  \begin{aligned}
    \mat R_b\vect v&=\vect v(\mathcal{I}_b),
      &&\vect v\in\mathbb{C}^{N_h},\;\vect v(\mathcal{I}_b)\in\mathbb{C}^{n_b},\\
    (\mat R_b^*\vect w)_j
      &=\begin{cases}
          w_\ell, &\text{if }j=i_\ell^{(b)}\text{ for some }\ell,\\
          0, &\text{otherwise},
        \end{cases}
      &&\vect w\in\mathbb{C}^{n_b},\;\mat R_b^*\vect w\in\mathbb{C}^{N_h}.
  \end{aligned}
  \end{equation*}
DD-POD constructs the local reduced bases via truncated SVDs of locally restricted snapshot matrices
\(\mat E_b=\mat R_b\mat E\approx\mat\Phi_b\mat\Sigma_b\mat V_b^*\), where \(\mat\Phi_b\in\mathbb C^{n_b\times r_b}\).
Prolonging the local bases to the full domain gives the reduced space
\begin{equation}
  \mathcal V^{\rm DD}
  =
  \bigoplus_{b=1}^B\range(\mat R_b^*\mat\Phi_b).
\label{eq:blockwise-space}
\end{equation}
The number of nonzero stored basis entries is
$\sum_{b=1}^B n_b r_b$.

\textbf{Subdomain partition.} In this paper,
the subdomains $\{\Omega_b\}_{b=1}^{B}$ are chosen a priori to represent
distinct wave behavior using independent local POD spaces. We align
the partition with known physical and geometric features, including
material interfaces, embedded structures, and the physical--PML
interface. Within the PML, subdomains are further subdivided by coordinate orientation and attenuation depth to separate directional propagation from the rapid exponential decay across the absorbing layer.  The specific
partitions used in the numerical experiments are described in
Section~\ref{sec:numerics}. Alternatively, data-driven partitioning through error-guided adaptive
refinement~\cite{carlberg2015adaptive,etter2020online} or clustering~\cite{rocha2020adaptive,ferreira2022adaptivity} can be used.

\textbf{Storage-aware greedy local rank selection.}
We select local ranks in each subdomain through a greedy procedure that aims to limit
basis storage while satisfying a prescribed global reconstruction
error tolerance.

Let $\sigma_{b,\ell}$ denote the singular values of $\mat E_b$.
Since the subdomains have disjoint DOF sets, the squared reconstruction
error is the sum of the squares of the locally discarded singular values:
\begin{equation}
  \|\mat E-\mat E^{\rm DD}\|_F^2
  =
  \sum_{b=1}^B\sum_{\ell>r_b}\sigma_{b,\ell}^2
  \leq \eta^2\|\mat E\|_F^2,
  \label{eq:block-weighted-error-bound}
\end{equation}
where $\mat E^{\rm DD}$ assembles the local truncated-SVD
approximations.

Adding the next mode in subdomain $b$ reduces the squared reconstruction
error by $\sigma_{b,r_b+1}^2$ at a basis storage cost of $n_b$ entries.
This motivates selecting the mode with the largest squared error
reduction per stored entry. 
Specifically, we initialize $r_b=0$ for
$1\leq b\leq B$ and, at each step, select
\begin{equation}
  b^*\in\operatorname*{arg\,max}_{b:\,r_b<\operatorname{rank}(\mat E_b)}
  \frac{\sigma_{b,r_b+1}^2}{n_b}.
  \label{eq:marginal-storage-greedy}
\end{equation}
We increase its rank, $r_{b^*}:= r_{b^*}+1$, and repeat this procedure
until \eqref{eq:block-weighted-error-bound} is satisfied.

\begin{remark}[Impact of domain decomposition on the total rank of the ROM]
\label{rem:spatial-localization}
The global POD constructs the reduced basis through truncated SVD, which gives the optimal linear reconstruction with respect to $||\cdot||_F$. Using local bases, we lose this optimality, and the total rank $\sum_{b=1}^Br_b\geq r_{\textrm{global}}$, where $r_{\textrm{global}}$ is the rank for the global POD. In practice, we observe that the rank growth is modest, while the memory reduction is substantial.
\end{remark}

\subsection{Online stage}
\label{sec:dd-rom-online}
Given the local bases \(\{\mat\Phi_b\}_{b=1}^B\) and local ranks
\(\{r_b\}_{b=1}^B\) produced by the offline stage, the online stage assembles
a Galerkin system of dimension \(r=\sum_b r_b\) and solves it for a new
parameter instance. Define the subdomain-prolonged basis matrix as
$$
\mat{\widetilde\Phi}= \begin{bmatrix}
		    \mat R_1^*\mat\Phi_{1}&\mat R_2^*\mat\Phi_{2}&\cdots&\mat R_B^*\mat\Phi_{B}
			  \end{bmatrix},\quad \mat{\widetilde\Phi}\in\mathbb{C}^{N_h\times r},\quad r=\sum_{b=1}^Br_b.
$$
The DD-POD Galerkin approximation seeks \(\vect u_{\mu}\approx \mat{\widetilde\Phi}\vect c_{\mu}\) by
solving the Galerkin system \(\mat{\widetilde\Phi}^* \mat A_{\mu}
\mat{\widetilde\Phi} \vect c_\mu=\mat{\widetilde\Phi}^*\vect b_\mu\).
Since each \(\mat\Phi_b\) is supported on \(\mathcal{I}_b\), the projection
operations decouple subdomain-wise:
\begin{subequations}
\begin{align}
&\mat{\widetilde\Phi}^*\vect b_\mu=\begin{pmatrix}\mat\Phi_1^*\mat R_1\vect b_{\mu},\dots,\mat\Phi_B^*\mat R_B\vect b_{\mu}
\end{pmatrix}^T
=\begin{pmatrix}\mat\Phi_1^*\vect b_{\mu}(\mathcal{I}_1),\dots,\mat\Phi_B^*\vect b_{\mu}(\mathcal{I}_B)
\end{pmatrix}^T,\\
&\mat{\widetilde\Phi}^*\mat A_{\mu}\mat{\widetilde\Phi}
  =
   \begin{bmatrix}
    \mat\Phi_1^* \mat A_{\mu}(\mathcal{I}_1,\mathcal{I}_1)\mat\Phi_1
    & \mat\Phi_1^* \mat A_{\mu}(\mathcal{I}_1,\mathcal{I}_2)\mat\Phi_2
    & \cdots
    & \mat\Phi_1^* \mat A_{\mu}(\mathcal{I}_1,\mathcal{I}_B)\mat\Phi_B
    \\[4pt]
    \mat\Phi_2^* \mat A_{\mu}(\mathcal{I}_2,\mathcal{I}_1)\mat\Phi_1
    & \mat\Phi_2^* \mat A_{\mu}(\mathcal{I}_2,\mathcal{I}_2)\mat\Phi_2
    & \cdots
    & \mat\Phi_2^* \mat A_{\mu}(\mathcal{I}_2,\mathcal{I}_B)\mat\Phi_B
    \\[4pt]
    \vdots & \vdots & \ddots & \vdots
    \\[4pt]
    \mat\Phi_B^* \mat A_{\mu}(\mathcal{I}_B,\mathcal{I}_1)\mat\Phi_1
    & \mat\Phi_B^* \mat A_{\mu}(\mathcal{I}_B,\mathcal{I}_2)\mat\Phi_2
    & \cdots
    & \mat\Phi_B^* \mat A_{\mu}(\mathcal{I}_B,\mathcal{I}_B)\mat\Phi_B
  \end{bmatrix}.
\label{eq:reduced-galerkin-block-matrix}
\end{align}
\end{subequations}
Partition the reduced coefficient vector as
\(\vect c_{\mu,r}=((\vect c_{\mu}^{(1)})^T,\dots,
(\vect c_{\mu}^{(B)})^T)^T\in\mathbb C^{r}\), where
\(\vect c_{\mu}^{(b)}\in\mathbb{C}^{r_b}\).  The reduced solution is then
  \[
  \vect u_{\mu}\approx\mat{\widetilde\Phi}\vect c_{\mu,r}
  =\sum_{b=1}^B\mat R_b^*\mat\Phi_b\vect c_{\mu}^{(b)}.
  \]

In this paper, we first assume that \(\mat A_\mu\) is affine, i.e., that it can be decomposed as
\(\mat A_\mu=\sum_{l=1}^{Q}\theta_l(\mu)\mat A_l\) with
\(\theta_l(\mu)\in\mathbb{C}\).  We can then precompute the projected operators
\(\mat A_{l,r}=\mat{\widetilde\Phi}^*\mat A_l\mat{\widetilde\Phi}
\in\mathbb{C}^{r\times r}\) following
\eqref{eq:reduced-galerkin-block-matrix}.  The online projected operator is
then assembled as
\(\mat{\widetilde\Phi}^*\mat A_\mu\mat{\widetilde\Phi}
=\sum_{l=1}^Q\theta_l(\mu)\mat A_{l,r}\).  The same offline--online
decomposition applies to an affine right-hand side. 
When the affine assumption does not hold, we directly perform projections online.
Further online acceleration for the nonaffine case can be achieved through hyper-reduction techniques such as the empirical interpolation method (EIM)~\cite{BARRAULT2004667} or discrete EIM (DEIM)~\cite{Chaturantabut2010DEIM}, which is left for future work.

In the online stage, subdomain-wise projection and prolongation require $O(\sum_{b=1}^B n_b r_b)$ operations, whereas assembling the $Q$-term affine reduced operator and factoring the dense reduced system scale as $O(Qr^2)$ and $O(r^3)$, respectively, with total reduced dimension $r=\sum_{b=1}^B r_b$. While decoupling the subdomains typically increases the aggregate dimension $r$ relative to a single global rank $r^{\rm G}$, the dimension of each local mode $n_b$ is much lower, $n_b \ll N_h$. Consequently, DD-POD reduces both memory footprint and projection cost whenever the basis compression $\sum_{b=1}^B n_b r_b \ll N_h r^{\rm G}$ outweighs the modest growth in the rank, as verified by numerical tests in Section~\ref{sec:numerics}.

\section{Trajectory-Aware DD-MSRB Preconditioner\label{sec:correction-rom}}
Preconditioners accelerate Krylov methods by providing inexpensive approximations to the inverse of the system operator. ROM-based preconditioners learn such approximations from parametric solution data~\cite{DalSanto2018,Peng2025,TangPeng2025}. In a flexible Krylov method, however, the preconditioned vector at each iteration determines the right-hand side of the subsequent inner system. On the one hand, each ROM is built from  inverse-action snapshots of the inner systems along the Krylov trajectory. On the other hand, the ROMs themselves determine the inner systems in subsequent iterations and thus the trajectory itself.  The MSRB approach accounts for this circular dependence between the iteration trajectory and the ROMs by constructing  ROMs sequentially, in a trajectory-aware manner~\cite{DalSanto2018,TangPeng2025}.  Here we adapt this trajectory-aware construction to
parametric Helmholtz--PML systems and represent each ROM with the
domain-decomposed bases of Section~\ref{sec:dd-rom}.

This section first reviews FGMRES, then defines the reduced inverse-action
preconditioner and its trajectory-aware offline construction.  Finally, the
constructed ROMs are combined with a shifted-Laplacian preconditioner to build  a
hybrid preconditioner for FGMRES. The offline and online procedures are
summarized in Algorithms~\ref{alg:offline-crom} and~\ref{alg:online-crom}.

\subsection{Flexible GMRES\label{sec:fgmres}}
To accommodate iteration-dependent reduced spaces, we adopt FGMRES~\cite{Saad1993}, which allows the right preconditioner $\mat M_j^{-1}$ to vary across iterations. FGMRES constructs the solution correction in the span of the preconditioned vectors $\vect z_j = \mat M_j^{-1}\vect v_j$. Suppressing explicit parameter dependence, we summarize the procedure in Algorithm~\ref{alg:flexible-gmres}, where $\langle\vect v,\vect w\rangle = \vect v^*\vect w$ denotes the standard Euclidean inner product on $\mathbb{C}^{N_h}$.

\begin{paperalgorithm}{Flexible GMRES}
{alg:flexible-gmres}
\textbf{Input:} system \(\mat A_\mu\vect u=\vect b_\mu\), initial state \(\vect u_0\), maximum depth
\(k\), right-preconditioning maps \(\mat M_j^{-1}\), \(j=1,\ldots,k\), and
FGMRES tolerance \(\tau\).

\textbf{Output:} FGMRES approximation \(\vect u_j\).

\begin{enumerate}
\item Form \(\vect r_0=\vect b_\mu-\mat A_\mu\vect u_0\).
      If \(\vect r_0=\vect 0\), return \(\vect u_0\).
      Otherwise define \(\vect q_1=\vect r_0/\|\vect r_0\|_2\).

\item For \(j=1,\ldots,k\):
      \begin{enumerate}

      \item Apply the current right preconditioner, \(\vect z_j=\mat M_j^{-1}\vect q_j.\)

      \item Compute \(\vect w^{(0)}=\mat A_\mu \vect z_j.\)

      \item For \(i=1,\ldots,j\), perform the Arnoldi orthogonalization:
            \[
              h_{ij}=\vect q_i^*\vect w^{(i-1)},
              \qquad
              \vect w^{(i)}=\vect w^{(i-1)}-h_{ij}\vect q_i.
            \]

      \item Set
            \(
              h_{j+1,j}=\|\vect w^{(j)}\|_2.
            \)
            If \(h_{j+1,j}=0\), mark Arnoldi breakdown and retain the
            current column for the least-squares update; otherwise define
            \(
              \vect q_{j+1}=\vect w^{(j)}/h_{j+1,j}.
            \)

      \item With \(\mat Z_j=[\vect z_1,\ldots,\vect z_j]\) and the current Hessenberg matrix
            \(\bar{\mat H}_j\), solve
            \[
              \vect y_j=\argmin_{\vect y}
              \bigl\|\|\vect r_0\|_2\vect e_1-\bar{\mat H}_j\vect y\bigr\|_2.
            \]
            Set \(\vect u_j=\vect u_0+\mat Z_j\vect y_j\) and evaluate
            \(\vect r_j=\vect b_\mu-\mat A_\mu\vect u_j\).
            Terminate if \(\|\vect r_j\|_2/\|\vect r_0\|_2\le\tau\) or
            breakdown was marked; report whether the true residual met
            the requested tolerance.

      \end{enumerate}

\item Return
      \(
        \vect u_j=\vect u_0+\mat Z_j\vect y_j,
      \)
      where \(j\le k\) is the final iteration index.

\end{enumerate}
\end{paperalgorithm}

\subsection{ROM-enhanced preconditioner and trajectory-aware ROM construction\label{sec:trajectory-aware-offline}}

Applying the exact inverse as a preconditioner yields immediate
convergence, but its application can be prohibitively expensive.
ROM-based preconditioners approximate the inverse action
$\mat A_\mu^{-1}\vect q_{\mu,j}$ in a reduced space. Given a parameter-independent reduced basis $\mat\Psi_j$ for
FGMRES step $j$, we define
\begin{equation}
  \widehat{\vect z}_{\mu,j}
  =\mat P_j(\mu)\vect q_{\mu,j}
  =\mat\Psi_j
  (\mat\Psi_j^*\mat A_\mu\mat\Psi_j)^{-1}
  \mat\Psi_j^*\vect q_{\mu,j}.
  \label{eq:rom-pc}
\end{equation}
The basis $\mat\Psi_j$ is shared across the parameter domain,
while $\mat P_j(\mu)$ depends on $\mu$ through the reduced operator.

A key observation is that the Arnoldi direction of the next iteration \(\vect q_{\mu,j+1}\) is
determined by the preconditioned vector of the current iteration \(\vect z_{\mu,j}\), and serves as
the right-hand side of the inner problem in the next iteration
\(\mat A_{\mu}\vect z_{\mu,j+1}=\vect q_{\mu,j+1}\)~\cite{DalSanto2018,TangPeng2025}.  Hence, the ROM for step \(j+1\) depends on
the preconditioners applied at steps \(1,\ldots,j\), so the spaces must be built
sequentially.  



For each $\mu\in\mathcal S_{\textrm{train}}$, we first compute a
reference solution $\vect u_\mu$ using shifted-Laplacian-preconditioned
GMRES or another suitable solver. For a prescribed number $k$ of
ROM-preconditioned iterations, we construct a sequence of iteration-dependent ROMs as follows.
\begin{enumerate}
\item \textbf{Initialize the first-iteration ROM-based preconditioner.}
Suppose the initial guess for FGMRES applied to \(\mat A_{\mu}\vect u_{\mu}=\vect b_{\mu}\) is \(\vect u_{\mu,0}\). 
For the first iteration, we have
\begin{equation}
  \vect q_{\mu,1}=\frac{\vect b_{\mu}-\mat A_{\mu} \vect u_{\mu,0}}{\|\vect b_{\mu}-\mat A_{\mu} \vect u_{\mu,0}\|_2},\; \textrm{then}\;
  \vect z_{\mu,1}=\mat A_\mu^{-1}\vect q_{\mu,1}=\frac{\vect u_{\mu}-\vect u_{\mu,0}}{\| \vect r_{\mu,0}\|_2}.
\end{equation}
The snapshots \(\vect z_{\mu,1}\), with
\(\mu\in\mathcal{S}_{\textrm{train}}\), are collected to form the training set.  Using these snapshots, we
construct the initial basis
\(\mat\Psi_1\) using DD-POD in Section~\ref{sec:dd-rom} and the associated preconditioner \(\mat P_1(\mu)\) defined in
\eqref{eq:rom-pc}.
For the zero initial guess used in this work,
\(\vect z_{\mu,1}=\vect u_\mu/\|\vect b_\mu\|_2\), so \(\mat\Psi_1\) is
built from scaled solution snapshots.

\item \textbf{Sequential ROM construction.}
For $j=2,\ldots,k$, we construct $\mat P_j$ by continuing the Arnoldi
process for each training parameter using the previously constructed ROM-based
preconditioners:
\begin{enumerate}

\item \textbf{Continue the Arnoldi process.} Apply the ROM preconditioner to the
offline direction,
\(\widehat{\vect z}_{\mu,j-1}=\mat P_{j-1}(\mu)\vect q_{\mu,j-1}\) defined in \eqref{eq:rom-pc}, and set
\(\vect w_{\mu,j-1}=\mat A_\mu\widehat{\vect z}_{\mu,j-1}\).  For
\(i=1,\ldots,j-1\), Gram--Schmidt orthogonalization gives
\begin{equation}
  \vect q_{\mu,j}=\frac{\vect w_{\mu,j-1}-\sum_{i=1}^{j-1}h_{i,j-1}^{(\mu)}\vect q_{\mu,i}}{h_{j,j-1}^{(\mu)}}.
\label{eq:trajectory-snapshot-recurrence-res}
\end{equation}
Given \(\vect z_{\mu,i}=\mat A_\mu^{-1}\vect q_{\mu,i}\) and \(\vect w_{\mu,j-1}=\mat A_{\mu}\widehat{\vect z}_{\mu,j-1}\),
applying \(\mat A_\mu^{-1}\) to both sides of \eqref{eq:trajectory-snapshot-recurrence-res} yields
\begin{align}
\vect z_{\mu,j}
  =\frac{\widehat{\vect z}_{\mu,j-1}-\sum_{i=1}^{j-1}h_{i,j-1}^{(\mu)}\vect z_{\mu,i}}
  {h_{j,j-1}^{(\mu)}}.
\label{eq:trajectory-snapshot-recurrence}
\end{align}
In exact arithmetic, every training parameter thus obtains the direction
\(\vect q_{\mu,j}\) and its inverse action \(\vect z_{\mu,j}\) without a new
full-order solve.
\item \textbf{Form \(\mat P_j\).}
 Collect the snapshot matrix \([\vect z_{\mu_1,j},\ldots,\vect z_{\mu_{n_s},j}]\), construct a DD-POD basis \(\mat\Psi_j\) following Section~\ref{sec:dd-rom}, and define \(\mat P_j\) by~\eqref{eq:rom-pc}.
\end{enumerate}

\end{enumerate}

\subsection{Hybrid preconditioner and offline--online algorithm}
\label{subsec:algorithm-summary}
To limit memory usage and offline overhead, ROMs are constructed only for the first $k$ FGMRES iterations, using a preselected set of training parameters. Beyond the $k$-th iteration, error modes not well resolved by the ROMs may stall convergence before the target tolerance is achieved. To address these error modes, we restart FGMRES after the ROM-preconditioned iterations and continue with the
shifted-Laplacian preconditioner.  The resulting hybrid preconditioner is
\begin{equation}
  \mat M_j^{-1}(\mu)=
  \begin{cases}
    \mat P_j(\mu)
    =\mat\Psi_j\bigl(\mat\Psi_j^*\mat A_\mu\mat\Psi_j\bigr)^{-1}\mat\Psi_j^*,
      & j=1,\ldots,\ell,\quad\text{first cycle (DD-MSRB)},\\[2mm]
    \mathcal B_{\mathrm{SL}}(\mu),
      & \text{all iterations after the restart (shifted-Laplacian)},
  \end{cases}
\label{eq:hybrid-pc}
\end{equation}
where \(\ell\le k\) is the number of ROM-preconditioned iterations taken, and
\(\mathcal B_{\mathrm{SL}}(\mu)\approx\mat M_\beta(\mu)^{-1}\) denotes
the shifted-Laplacian preconditioner \(\mathcal B_\beta\) of
Section~\ref{subsec:helmholtz-solvers} at the parameter \(\mu\).  We call the
iterations of the first cycle ROM steps and the iterations after the restart
shifted-Laplacian (SL) steps.  The ROM steps reduce the dominant
residual components by leveraging low-dimensional structures across parameters, and the
SL steps reduce the remainder.

As in Algorithm~\ref{alg:flexible-gmres}, the true residual of the assembled
PML system is evaluated after every ROM step, and the solve terminates as soon
as it meets the tolerance~\(\tau\).  Otherwise, the ROM steps continue until
all \(k\) reduced spaces have been used or the Arnoldi process breaks down,
which defines \(\ell\le k\).  FGMRES is then restarted from the ROM iterate
\(\vect u_\ell\).  Because \(\mathcal B_{\mathrm{SL}}(\mu)\) does not change
across iterations, after restarting, the SL steps coincide with right-preconditioned
GMRES(\(m\)) started from \(\vect u_\ell\).  The shifted-Laplacian
preconditioner is set up only at this restart; a query that converges within
the ROM steps therefore incurs no shifted-Laplacian setup cost.
We summarize the offline training to construct the DD-MSRB preconditioner in Algorithm~\ref{alg:offline-crom} and the online solve with the hybrid preconditioner in
Algorithm~\ref{alg:online-crom}.

\begin{remark}
We combine the DD-MSRB preconditioner with a shifted-Laplacian preconditioner.
Other classical preconditioners can be used instead, as long as they
effectively address the residual components not well represented by the ROMs.
\end{remark}

\begin{paperalgorithm}{Offline trajectory-aware ROM training}
{alg:offline-crom}
\textbf{Input:} training parameters \(\{\mu_i\}_{i=1}^{n_s}\), systems
\(\mat A_i=\mat A_{\mu_i}\), right-hand sides \(\vect b_i\), initial guesses \(\vect u_{0,i}\),
reference solutions \(\vect u_i^\star\) from a verified full-order solver,
number of ROM-preconditioned iterations \(k\), SVD truncation rule.

\textbf{Output:} bases \(\mat\Psi_1,\ldots,\mat\Psi_k\) defining the ROM
preconditioners \(\mat P_1,\ldots,\mat P_k\) in~\eqref{eq:rom-pc}.

\begin{enumerate}
\item For each training parameter, form the initial residual
      \(\vect r_{0,i}=\vect b_i-\mat A_i \vect u_{0,i}\) and set \(\vect q_{1,i}=\vect r_{0,i}/\|\vect r_{0,i}\|_2\).
\item Form the first snapshots
      \(\vect z_{1,i}=(\vect u_i^\star-\vect u_{0,i})/\|\vect r_{0,i}\|_2\)
      and compress \([\vect z_{1,1},\ldots,\vect z_{1,n_s}]\) into \(\mat P_1\).
\item For \(j=1,\ldots,k-1\):
      \begin{enumerate}
      \item For every training parameter \(\mu_i\), \(i=1,\ldots,n_s\), apply
            \(\widehat{\vect z}_{j,i}=\mat P_j(\mu_i)\vect q_{j,i}\), set
            \(\vect w_{j,i}=\mat A_i\widehat{\vect z}_{j,i}\), and
            orthogonalize against
            \(\vect q_{1,i},\ldots,\vect q_{j,i}\) to obtain the Hessenberg
            entries \(h_{1j}^{(i)},\ldots,h_{j+1,j}^{(i)}\) and
            \(\vect q_{j+1,i}\).
      \item For every training parameter without numerical Arnoldi breakdown, infer
            \(\vect z_{j+1,i}=\mat A_i^{-1}\vect q_{j+1,i}\) from
            \eqref{eq:trajectory-snapshot-recurrence}, verify its residual,
            and replace it by a verified full-order solve if necessary.
      \item After all training parameters have been advanced, compress
            \([\vect z_{j+1,1},\ldots,\vect z_{j+1,n_s}]\) once to construct
            \(\mat\Psi_{j+1}\) and \(\mat P_{j+1}\).
      \end{enumerate}
\item Store the bases, local ranks, and reduced assembly data.
\end{enumerate}
\end{paperalgorithm}

\begin{paperalgorithm}{Online FGMRES with the hybrid preconditioner}
{alg:online-crom}
\textbf{Input:} parameter \(\mu\), initial guess \(\vect u_0\), trained bases
\(\{\mat\Psi_j\}_{j=1}^k\), restart length \(m\), and tolerance~\(\tau\).

\textbf{Output:} approximation of \(\vect u(\mu)\) with verified true residual.

\begin{enumerate}
\item Assemble \(\mat A_{\mu}\) and \(\vect b_\mu\). Form
      \(\vect r_0=\vect b_\mu-\mat A_{\mu}\vect u_0\).
      If \(\vect r_0=\vect 0\), return \(\vect u_0\);
      otherwise set \(\vect q_1=\vect r_0/\|\vect r_0\|_2\).
\item \textbf{ROM steps.} For \(j=1,\ldots,k\):
      \begin{enumerate}
      \item Precondition with the reduced Galerkin solve
            in~\eqref{eq:rom-pc}: \(\vect z_j=\mat P_j(\mu)\vect q_j\).
      \item Form \(\vect w_j=\mat A_{\mu}\vect z_j\), apply two-pass
            Gram--Schmidt orthogonalization against \(\vect q_1,\ldots,\vect q_j\),
            and normalize to obtain \(\vect q_{j+1}\).
      \item Solve
            \(\vect y_j=\argmin_{\vect y}
            \bigl\|\|\vect r_0\|_2\vect e_1-\bar{\mat H}_j\vect y\bigr\|_2\),
            set \(\vect u_j=\vect u_0+[\vect z_1,\ldots,\vect z_j]\vect y_j\), and
            evaluate the true residual
            \(\vect r_j=\vect b_\mu-\mat A_{\mu}\vect u_j\).
            If \(\|\vect r_j\|_2/\|\vect r_0\|_2\le\tau\), return \(\vect u_j\).
            If numerical breakdown was detected, leave the loop.
      \end{enumerate}
\item \textbf{SL steps.} With \(\ell\le k\) the number of ROM steps taken,
      set up \(\mathcal B_{\mathrm{SL}}(\mu)\), restart from \(\vect u_\ell\),
      and apply right-preconditioned GMRES(\(m\)) with
      \(\mat M^{-1}(\mu)=\mathcal B_{\mathrm{SL}}(\mu)\)
      until the true residual meets \(\tau\).
\item Return the final iterate with its verified true residual.
\end{enumerate}
\end{paperalgorithm}

\section{Numerical Experiments}\label{sec:numerics}
We assess the hybrid preconditioner~\eqref{eq:hybrid-pc} in three experiments.
The two-dimensional half-space problem establishes the basic effects of
trajectory-aware training and of domain decomposition in DD-MSRB.  Elliptic-obstacle
scattering tests a problem with parametrized geometry.  The
three-dimensional Gaussian waveguide examines storage and online efficiency at
large scale.

Each experiment compares DD-MSRB with the shifted-Laplacian baseline and with
global MSRB, which uses a single global POD basis per iteration and the same
trajectory-aware construction.  Unless
stated otherwise, \(k=3\) ROM-preconditioned iterations are used; global MSRB and DD-MSRB use
the same training and test parameters and iteration-wise truncation tolerances
\(\eta_j\).  All offline and online solves start from the zero initial guess;
when the ROM steps do not reach the tolerance, FGMRES restarts from the ROM
iterate \(\vect u_\ell\) for the SL steps, and the SL baseline starts from zero.
We report the total online time, which excludes
full-order assembly and offline training, and ROM storage as both the payload
and the size of the deployed MATLAB objects.  The payload counts the stored
reduced bases, precomputed reduced operators, and, for DD, one shared subdomain
index map, as in Section~\ref{sec:dd-rom-online}.  We also report the true
relative residual over the full PML-extended domain,
\begin{equation}
  \rho
  =
  \frac{\|\vect b_\mu-\mat A_{\mu}\vect u\|_2}{\|\vect b_\mu\|_2}.
\label{eq:true-residual}
\end{equation}  We denote the relative residual at the end of
the ROM steps, before any SL steps, by \(\rho_{\rm ROM}\).

In the two-dimensional experiments, sparse direct solves provide the training
references and a second baseline, and \(\mathcal B_{\mathrm{SL}}\) is an ILU
factorization of the SL operator with the shift parameter \(\beta=0.2\)  and ILU
drop tolerance \(10^{-3}\), except in the training-coverage study
(Table~\ref{tab:jsc-halfspace-reuse}), which uses an exact sparse LU
factorization of the SL operator.  The SL steps use restart length \(m=40\), at most
20 cycles, and two-pass orthogonalization.  The SL baseline is labeled Shifted
ILU-GMRES in the tables and figures.  These computations use MATLAB
R2025a on an Apple M4 Max (16 cores, 128 GB) with a single thread; unless
stated otherwise, each query is timed five times after one untimed warm-up
run, with the method order rotated between repetitions so that no method is
always run first.  Offline times are the median wall-clock time of three
training runs after an untimed warm-up, including the construction of the DD
partition.  

In the three-dimensional experiments, the shifted-Laplacian preconditioner uses the shift
\(\beta=2\) and a six-level multigrid V-cycle with two pre- and post-smoothing
Jacobi sweeps (weight 0.6) and at most 12,000 coarse unknowns; the SL steps use
restart length \(m=40\) and at most 40 cycles.  The training inverse actions are
computed with shifted-Laplacian GMRES to a relative residual of \(10^{-7}\); the
training solves of each iteration are independent and run in parallel on
four-thread MATLAB workers.  Online queries are solved one
at a time with at most eight threads, all meet \(\rho\le10^{-6}\), and each is
timed once unless stated otherwise.  These computations use MATLAB R2024a on
Hygon 7285 and 7285H CPU nodes with 64 cores and up to 224 GiB of memory per job.

\subsection{Experiment 1: Half-space Helmholtz problem}
\label{sec:helm-correction-rom}
We adapt the two-material cavity of~\cite{ChenHesthavenMadayRodriguez2010} to the scalar
Helmholtz--PML system of Section~\ref{sec:problem}.  We first examine the effects of trajectory-aware training
and the DD representation, then assess sensitivity to the
training set size, spatial partition, and frequency.

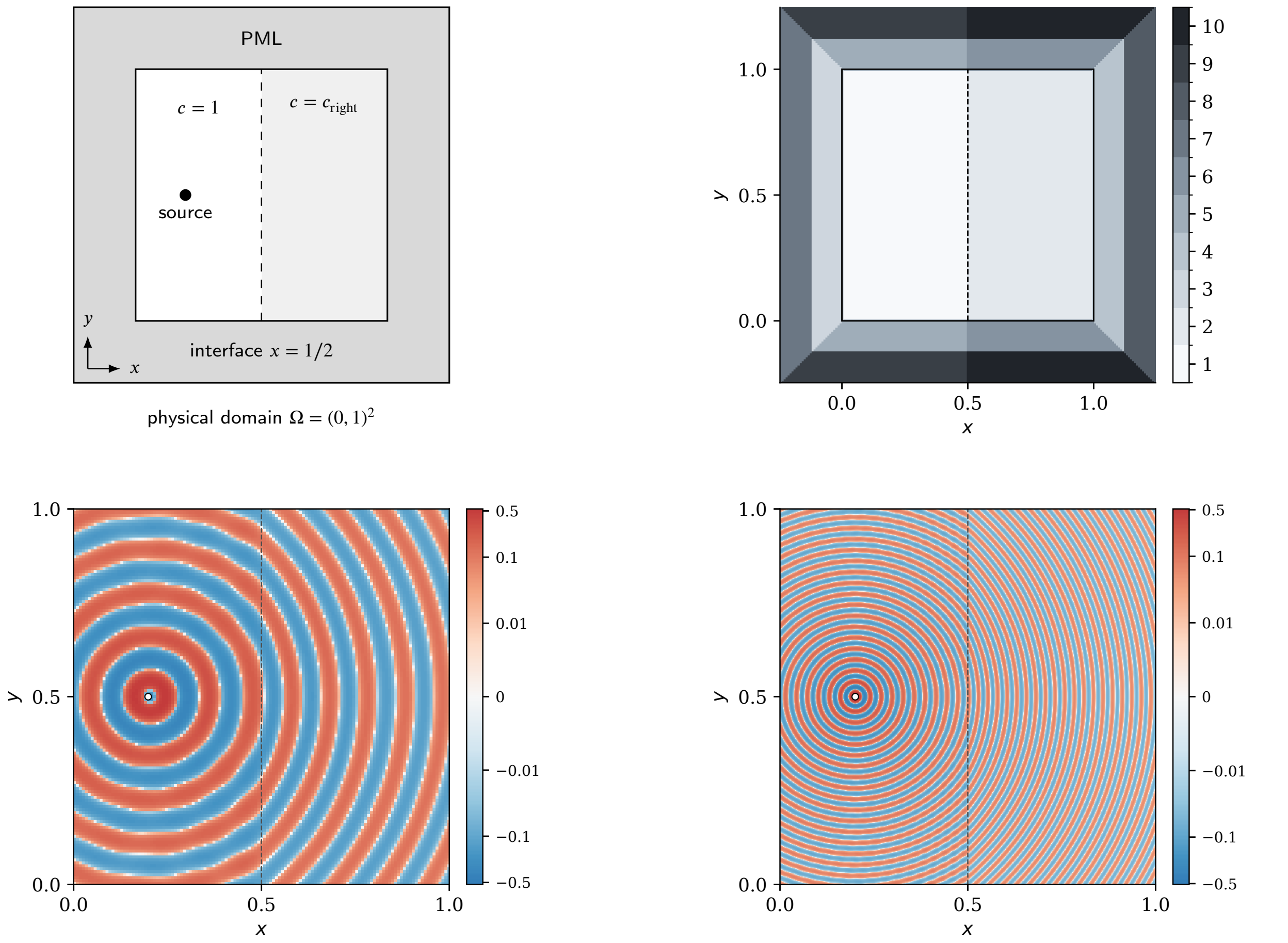
\begin{figure}[t]
\centering
\setlength{\panelunit}{0.1022727\linewidth}
\begin{tikzpicture}[
  x=\panelunit, y=\panelunit,
  src/.style={circle, fill=black, inner sep=1.5pt},
  arr/.style={-{Latex[length=1.5mm]}, line width=0.5pt},
  lab/.style={font=\scriptsize, align=center}
]
\useasboundingbox (0,0) rectangle (4.4,3.8);
\begin{scope}[shift={(1.0436,1.0036)}, scale=1.9168]
\fill[gray!30] (-0.246,-0.246) rectangle (1.246,1.246);
\fill[white] (0,0) rectangle (0.5,1);
\fill[gray!12] (0.5,0) rectangle (1,1);
\draw[line width=0.6pt] (0,0) rectangle (1,1);
\draw[dashed, line width=0.5pt] (0.5,0) -- (0.5,1);
\draw[line width=0.6pt] (-0.246,-0.246) rectangle (1.246,1.246);

\node[src] at (0.198,0.5) {};
\node[lab, below=2pt] at (0.198,0.5) {source};

\node[lab] at (0.25,0.85) {\(c=1\)};
\node[lab] at (0.75,0.85) {\(c=c_{\rm right}\)};
\node[lab] at (0.5,-0.123) {interface \(x=1/2\)};
\node[lab] at (0.5,1.123) {PML};

\draw[arr] (-0.19,-0.19) -- (-0.06,-0.19) node[lab, right] {\(x\)};
\draw[arr] (-0.19,-0.19) -- (-0.19,-0.06) node[lab, above] {\(y\)};
\end{scope}
\node[lab] at (2.002,0.27) {physical domain \(\Omega=(0,1)^2\)};
\end{tikzpicture}%
\hfill
\includegraphics[width=0.45\linewidth]{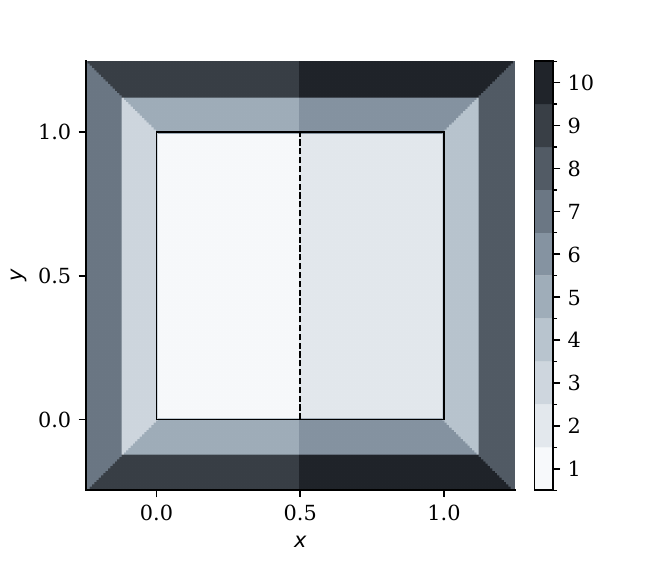}\par
\includegraphics[width=0.45\linewidth]{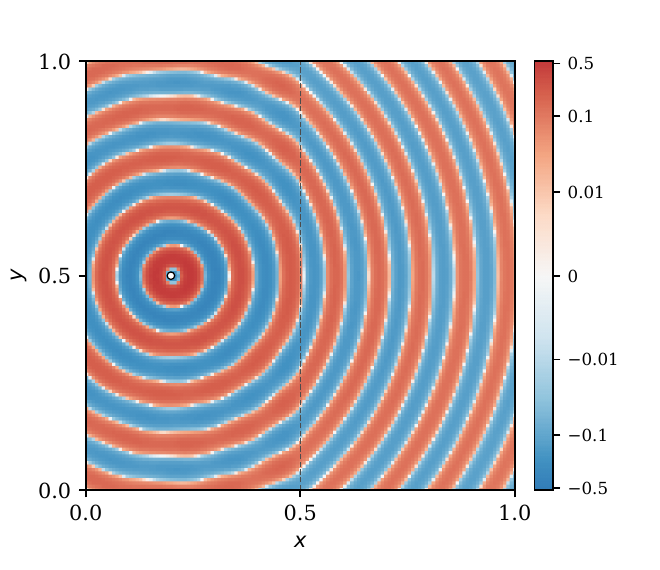}%
\hfill
\includegraphics[width=0.45\linewidth]{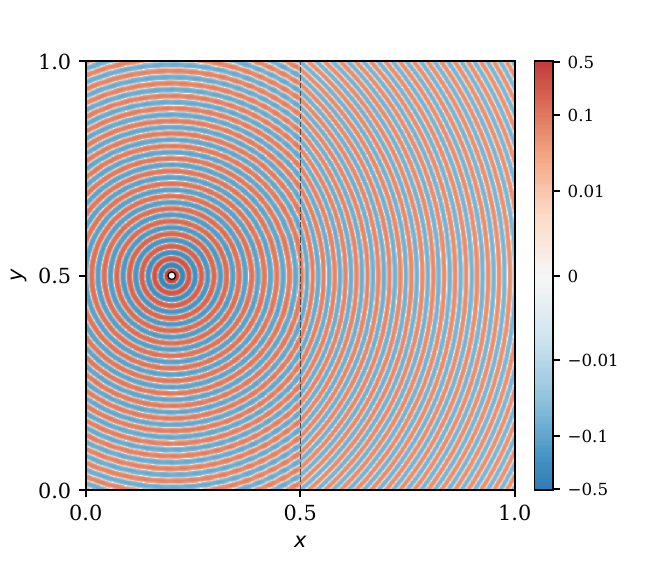}
\caption{Half-space problem in Section~\ref{sec:helm-correction-rom}.  Top left: model
geometry.  Top right: ten-subdomain DD partition.  Bottom left: normalized
real part of the solution in the physical domain at \(\omega=53.5\) and
\(c_{\rm right}=0.8622\).  Bottom right: the same at \(\omega=214\) and
\(c_{\rm right}=0.8761\).}
\label{fig:halfspace-geometry}
\end{figure}

\experimentlabel{Problem setting}
The physical domain is \(\Omega=(0,1)^2\), and the wavespeed is piecewise constant,
\begin{equation}
  c(x,y)=
  \begin{cases}
    1, & x<\tfrac{1}{2},\\
    c_{\rm right}, & x\geq\tfrac{1}{2},
  \end{cases}
\label{eq:halfspace-wavespeed}
\end{equation}
with parameter \(\mu=c_{\rm right}\).  A
point source is located at \((x_s,y_s)\approx(0.1984,0.50)\) in the left half.
The physical grid has 125 points per side, \(h=1/126\), and the PML has 32
cells on each side with \(\sigma_{\max}=80\), which gives a \(189\times189\)
grid and \(N_h=35{,}721\). Figure~\ref{fig:halfspace-geometry} shows the material configuration (top left) and the normalized real part of the solution at \(\omega=53.5\) and \(c_{\rm right}=0.8622\) (bottom left).  Unless otherwise specified, the angular frequency is
\(\omega=53.5\).

\begin{figure}
\centering
\includegraphics[width=.49\linewidth]{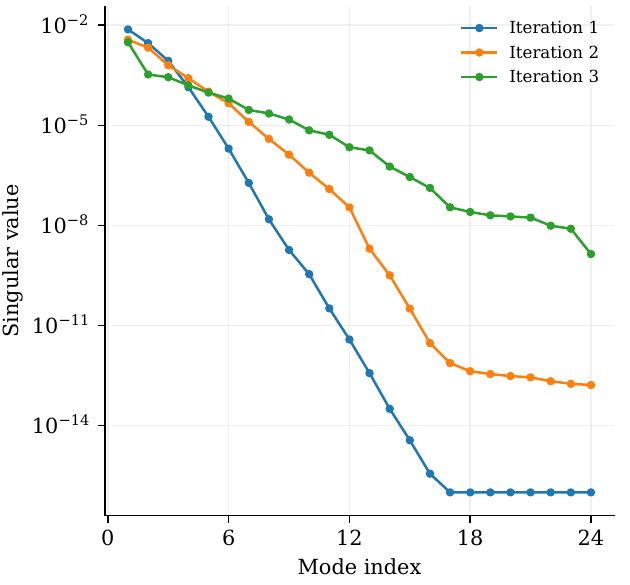}
\includegraphics[width=.49\linewidth]{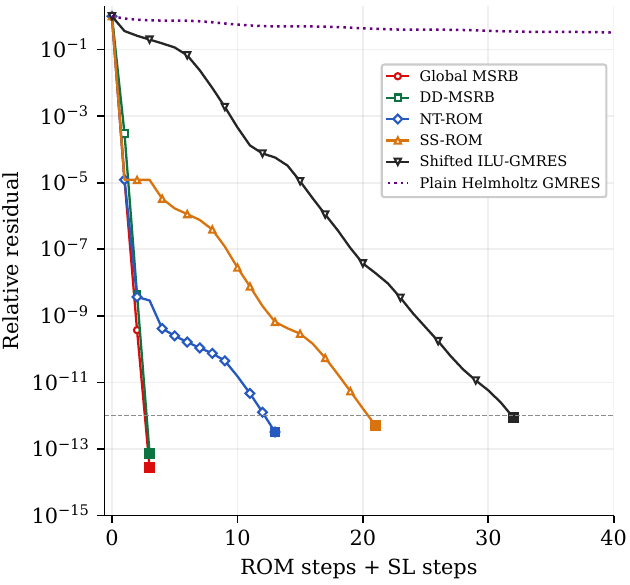}
\caption{Half-space problem in Section~\ref{sec:helm-correction-rom}.  Left: singular
values of the iteration-wise global snapshots.  Right: relative residual
histories at \(c_{\rm right}=0.8622\).}
\label{fig:helm-singular-values}
\end{figure}

Training uses \(n_s=24\) values of \(c_{\rm right}\) uniformly spaced in
\([0.80,0.95]\), and the test set contains 23 randomly sampled parameters, one in each training
interval; the representative query is \(c_{\rm right}=0.8622\).  Global MSRB and DD-MSRB both
use \(\eta_j=10^{-5}\) to construct ROMs for the first three iterations, guided by the snapshot spectra in
Figure~\ref{fig:helm-singular-values}.  The relative residual tolerance is \(\rho\le10^{-12}\) for all methods.  To
assess trajectory-aware training, we consider two global ROMs built from
the same training set without using the trajectory-aware construction approach. A non-trajectory-aware ROM (NT-ROM) is built iteration by iteration but uses the snapshots provided by
the unpreconditioned Arnoldi process, while a solution-snapshot ROM
(SS-ROM) uses full-order solution snapshots to construct a
single Galerkin space reused at every ROM-preconditioned step. Both are applied in the same hybrid preconditioner framework as the trajectory-aware reduced spaces of MSRB to isolate the effect of trajectory awareness.

By default, we use the ten-subdomain DD partition in Figure~\ref{fig:halfspace-geometry} (top right), which separates the two material regions (subdomains 1--2) and divides the PML by orientation and attenuation depth (inner subdomains 3--6 and outer subdomains 7--10), allowing local ranks to adapt to spatially heterogeneous wave behavior. The partition is shared by all training and test parameters. The sensitivity to the partition is also tested, and the alternative partitions are presented with that study.

\begin{table}[!htbp]\centering
\centering\normalsize\setlength{\tabcolsep}{2.5pt}
\setlength{\belowcaptionskip}{6pt}
\caption{Half-space problem in Section~\ref{sec:helm-correction-rom}, online cost (s)
over the 23 test queries.  It.\ is the mean number of SL steps per query
(GMRES iterations for Shifted ILU-GMRES); Speedup is the median time of
Shifted ILU-GMRES divided by that of each method.}
\label{tab:helm-cost}
\begin{tabular}{@{}lrrrr@{}}\toprule
Method & Mean & Median & It. & Speedup \\
\midrule
Global MSRB & 0.0072 & 0.0070 & 0.00 & $59.7\times$ \\
DD-MSRB & 0.0273 & 0.0078 & 0.09 & $53.6\times$ \\
NT-ROM & 0.2853 & 0.2880 & 9.57 & $1.5\times$ \\
SS-ROM & 0.3773 & 0.3890 & 17.57 & $1.1\times$ \\
Shifted ILU-GMRES & 0.4500 & 0.4177 & 32.13 & $1.0\times$ \\
Direct & 0.0551 & 0.0551 & -- & $7.6\times$ \\
\bottomrule\end{tabular}
\end{table}

\experimentlabel{Effect of trajectory-aware training}
Trajectory-aware training substantially reduces the iterations needed beyond
the ROM steps.  For all 23 test parameters, global MSRB reaches the tolerance
\(\rho\le10^{-12}\) within the three ROM steps. Trajectory-aware DD-MSRB converges within the first three ROM-preconditioned iterations for
21 tests, needing only one SL step for the other two
(Table~\ref{tab:helm-cost}).  NT-ROM and SS-ROM instead require an average of
\(9.57\) and \(17.57\) SL steps.  At the representative parameter
\(c_{\rm right}=0.8622\), NT-ROM, SS-ROM, and global MSRB share the first
reduced space for the zero initial guess, so their residuals coincide after the
first step (Figure~\ref{fig:helm-singular-values}).  The histories separate
afterwards: both MSRB variants reach the tolerance within three ROM steps, while
NT-ROM and SS-ROM continue with SL steps.  The SL baseline needs 32 iterations, and the relative residual of unpreconditioned GMRES is still above $10^{-1}$ after 40 iterations. 

The saved SL steps translate into online speedups.  Global MSRB and DD-MSRB
achieve median speedups of \(59.7\times\) and \(53.6\times\) over the SL
baseline, compared with \(1.5\times\) for NT-ROM and \(1.1\times\) for SS-ROM
(Table~\ref{tab:helm-cost}).  These large speedups come mainly from skipping the
SL stage: a test that converges within the ROM steps needs neither the SL
setup nor any SL steps, which take median times of \(0.2171\) and \(0.2005\) s,
respectively, out of the \(0.4177\) s of the SL baseline.  NT-ROM and SS-ROM still pay the
SL setup, which alone limits their speedups to about \(1.9\times\).  The median
online times of global MSRB and DD-MSRB, \(0.0070\) and \(0.0078\) s, are seven
to eight times shorter than that of the sparse direct solver.

\begin{figure}
\centering
\includegraphics[width=0.45\linewidth]
{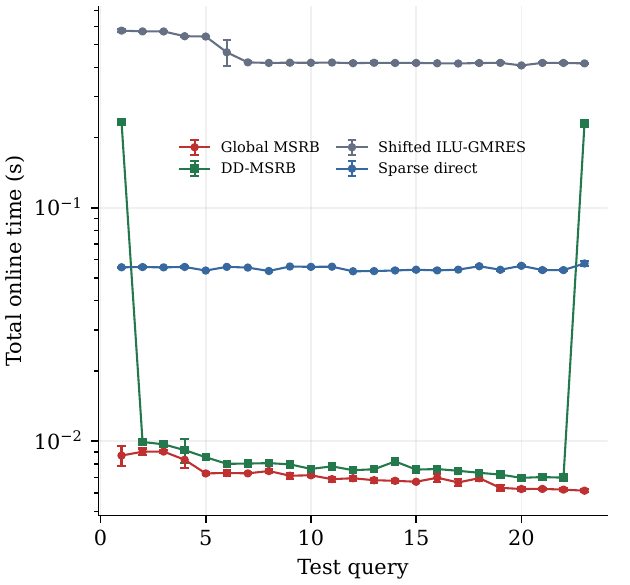}
\hfill
\includegraphics[width=0.45\linewidth]
{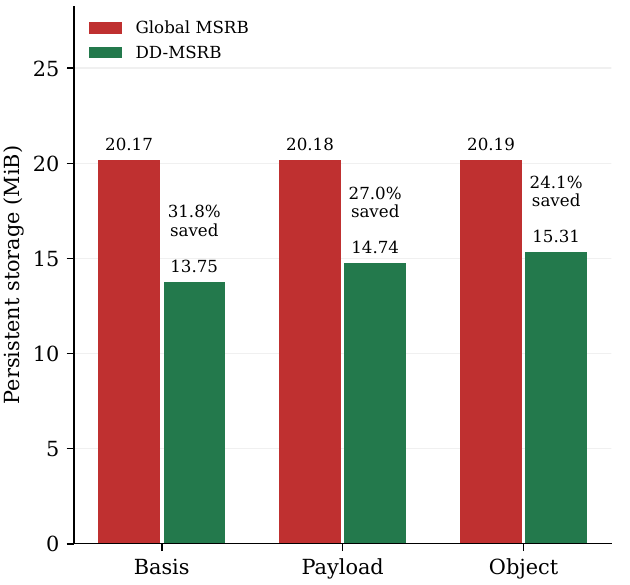}
\caption{Half-space problem in Section~\ref{sec:helm-correction-rom}, computational
costs.  Left: online times over the 23 test queries (mean \(\pm\) standard deviation over five
repetitions); DD queries 1 and 23 require the SL preconditioner setup and one
SL step each.  Right: total basis storage, deployed payload, and MATLAB object
size.}
\label{fig:halfspace-cost}
\end{figure}

\experimentlabel{Storage and online cost of DD}
DD-MSRB reduces the storage at an online cost comparable to global MSRB.  Under the same
truncation tolerances, it reduces the total basis storage from \(20.167\) to
\(13.750\) MiB, a saving of \(31.8\%\) (Figure~\ref{fig:halfspace-cost}, right).  The
deployed payload and the MATLAB object size shrink by \(27.0\%\) and
\(24.1\%\), slightly less than the bases alone, because DD stores larger dense reduced
operators and a subdomain index map.  The median online time of DD-MSRB stays
close to that of global MSRB (Figure~\ref{fig:halfspace-cost}, left); its mean
is higher because the end-point queries 1 and 23 need an SL step and thus pay
the SL preconditioner setup.

\experimentlabel{Sensitivity to the training set size}
Table~\ref{tab:jsc-halfspace-reuse} compares ROMs trained on \(n_s=6\), 12,
and 24 parameters in the same interval and benchmarked on nine fixed but randomly picked test parameters.

The storage advantage of DD also depends on the number of training samples.
At \(n_s=6\), both global MSRB and DD-MSRB retain all snapshot modes, leaving
no basis-storage savings.  The independent local coefficients increase the
aggregate reduced dimension and hence the storage of the dense reduced
operators, contributing to a \(40.0\%\) increase in deployed payload.  With
\(n_s=12\) and \(24\), the parameter range is better resolved, and DD-MSRB
instead reduces the payload by \(13.0\%\) and \(27.0\%\), respectively.

As expected, increasing the training set size \(n_s\) reduces the maximum
residual after the ROM steps and the mean online solve time
(Table~\ref{tab:jsc-halfspace-reuse}).  For the same training set size,
DD-MSRB and global MSRB have comparable online costs, while DD-MSRB requires
less offline construction time with sufficiently large training sets.

\begin{table}[t]\centering\normalsize
\setlength{\belowcaptionskip}{6pt}
\caption{Half-space problem in Section~\ref{sec:helm-correction-rom}, sensitivity to
the training set size \(n_s\) over nine fixed queries, two of them outside the
training interval.}\label{tab:jsc-halfspace-reuse}
\setlength{\tabcolsep}{3pt}
\begin{tabular}{rlrrrrr}\toprule
$n_s$ & Method & Payload (MiB) & Saving & Offline (s) & Time (s) & Max. $\rho_{\rm ROM}$ \\
\midrule
6 & Global MSRB & 6.543 & \(0.0\%\) & 0.74 & 0.2270 & \(7.14\times10^{-3}\) \\
6 & DD-MSRB & 9.162 & \(-40.0\%\) & 0.76 & 0.2679 & \(5.84\times10^{-3}\) \\
12 & Global MSRB & 15.816 & \(0.0\%\) & 2.12 & 0.1373 & \(2.15\times10^{-6}\) \\
12 & DD-MSRB & 13.766 & \(13.0\%\) & 1.54 & 0.1345 & \(1.82\times10^{-6}\) \\
24 & Global MSRB & 20.183 & \(0.0\%\) & 4.28 & 0.0359 & \(2.14\times10^{-9}\) \\
24 & DD-MSRB & 14.739 & \(27.0\%\) & 2.96 & 0.0454 & \(7.98\times10^{-9}\) \\
\bottomrule\end{tabular}\end{table}

\experimentlabel{Sensitivity to the spatial partition}
Adapting the spatial partition to material interfaces and PML geometry yields
substantially larger storage savings than naive or uniform subdivisions.
Table~\ref{tab:reviewer-partition} compares four partitions, shown in
Figure~\ref{fig:halfspace-partitions}:
(i) a three-subdomain (merged-PML) partition consisting of the two material
regions and the entire PML;
(ii) a ten-subdomain signed/depth partition consisting of the two material
regions and eight PML subdomains obtained by splitting each of the four PML
sides into inner and outer depth layers;
(iii) the ten-subdomain material/depth partition in
Figure~\ref{fig:halfspace-geometry} (top right), which aligns the top and
bottom PML subdivisions with the physical material interface; and
(iv) a uniform \(4\times4\) partition of the entire extended domain.
All four DD models are trained on the same 24 parameters with
\(\eta_j=10^{-5}\), and the corresponding DD-MSRB preconditioners reach the
prescribed \(10^{-10}\) tolerance within three ROM-preconditioned steps for all
five test queries with comparable online times.  However, their storage
savings differ.

Decomposing the PML by orientation and depth increases the payload saving from
\(11.5\%\) for the merged-PML partition to \(22.5\%\) for the signed/depth
partition. With the same ten-subdomain count (Figure~\ref{fig:halfspace-partitions}), the material/depth partition further increases the savings to \(27.0\%\) (versus \(22.5\%\) for the signed/depth partition) by splitting the top and bottom PMLs along the material interface, so that no PML subdomain straddles both media.   Conversely, increasing the subdomain count
without physical alignment can degrade compression.  The 16-subdomain uniform
\(4\times4\) tiling reduces the payload by only \(19.8\%\), because its tiles
straddle the boundary between the physical domain and the PML, which inflates
the aggregate reduced dimension and, consequently, the storage overhead of the
dense reduced operators.

\begin{table}[t]\centering\normalsize
\setlength{\belowcaptionskip}{6pt}
\caption{Half-space problem in Section~\ref{sec:helm-correction-rom}, sensitivity to
the spatial partition over five fixed queries.}\label{tab:reviewer-partition}
\setlength{\tabcolsep}{3pt}
\begin{tabular}{llrrrr}\toprule
Partition & Iteration-wise ranks & Payload (MiB) & Saving & Time (s) & Max. $\rho_{\rm ROM}$ \\
\midrule
Global & \((7,11,19)\) & 20.183 & 0.0\% & 0.0057 & \(1.70\times10^{-13}\) \\
Material/depth & \((44,80,123)\) & 14.739 & 27.0\% & 0.0069 & \(1.91\times10^{-12}\) \\
Merged PML & \((17,29,44)\) & 17.854 & 11.5\% & 0.0066 & \(2.16\times10^{-13}\) \\
Signed/depth & \((49,87,132)\) & 15.641 & 22.5\% & 0.0060 & \(1.43\times10^{-12}\) \\
Uniform tiles & \((72,137,201)\) & 16.185 & 19.8\% & 0.0066 & \(9.20\times10^{-13}\) \\
\bottomrule\end{tabular}

\bigskip
\tikzset{
  hspart/.style={scale=0.95, font=\tiny},
  hspml/.style={fill=gray!30},
  hsleft/.style={fill=white},
  hsright/.style={fill=gray!12},
  hscut/.style={line width=0.4pt},
}
\newcommand{\hsbase}{%
  \fill[hspml] (-0.5,-0.5) rectangle (2.5,2.5);
  \fill[hsleft] (0,0) rectangle (1,2);
  \fill[hsright] (1,0) rectangle (2,2);
  \draw[line width=0.7pt] (-0.5,-0.5) rectangle (2.5,2.5);
  \draw[line width=0.7pt] (0,0) rectangle (2,2);
  \draw[dashed, line width=0.6pt] (1,0) -- (1,2);
  \node at (0.5,1) {1};
  \node at (1.5,1) {2};}
\newcommand{\hsdepth}{%
  \draw[hscut] (-0.25,-0.25) rectangle (2.25,2.25);
  \draw[hscut] (0,0) -- (-0.5,-0.5) (2,0) -- (2.5,-0.5)
               (0,2) -- (-0.5,2.5) (2,2) -- (2.5,2.5);}
\begin{tikzpicture}[hspart]
\hsbase \hsdepth
\draw[hscut] (1,2) -- (1,2.5) (1,0) -- (1,-0.5);
\node at (-0.125,1) {3}; \node at (2.125,1) {4};
\node at (0.5,2.125) {5}; \node at (0.5,-0.125) {5};
\node at (1.5,2.125) {6}; \node at (1.5,-0.125) {6};
\node at (-0.375,1) {7}; \node at (2.375,1) {8};
\node at (0.5,2.375) {9}; \node at (0.5,-0.375) {9};
\node at (1.5,2.375) {10}; \node at (1.5,-0.375) {10};
\node[font=\footnotesize, below] at (1,-0.55) {(a) Material/depth};
\end{tikzpicture}\hfill
\begin{tikzpicture}[hspart]
\hsbase
\node at (1,2.25) {3};
\node[font=\footnotesize, below] at (1,-0.55) {(b) Merged PML};
\end{tikzpicture}\hfill
\begin{tikzpicture}[hspart]
\hsbase \hsdepth
\node at (-0.125,1) {3}; \node at (2.125,1) {4};
\node at (1,-0.125) {5}; \node at (1,2.125) {6};
\node at (-0.375,1) {7}; \node at (2.375,1) {8};
\node at (1,-0.375) {9}; \node at (1,2.375) {10};
\node[font=\footnotesize, below] at (1,-0.55) {(c) Signed/depth};
\end{tikzpicture}\hfill
\begin{tikzpicture}[hspart]
\fill[hspml] (-0.5,-0.5) rectangle (2.5,2.5);
\fill[hsleft] (0,0) rectangle (1,2);
\fill[hsright] (1,0) rectangle (2,2);
\draw[line width=0.7pt] (-0.5,-0.5) rectangle (2.5,2.5);
\draw[line width=0.7pt] (0,0) rectangle (2,2);
\foreach \c in {0.25,1,1.75} {
  \draw[hscut] (\c,-0.5) -- (\c,2.5);
  \draw[hscut] (-0.5,\c) -- (2.5,\c);}
\foreach \i [count=\ix from 0] in {-0.125,0.625,1.375,2.125} {
  \foreach \j [count=\jy from 0] in {-0.125,0.625,1.375,2.125} {
    \pgfmathtruncatemacro{\k}{1+\ix+4*\jy}
    \node at (\i,\j) {\k};}}
\node[font=\footnotesize, below] at (1,-0.55) {(d) Uniform \(4\times4\)};
\end{tikzpicture}

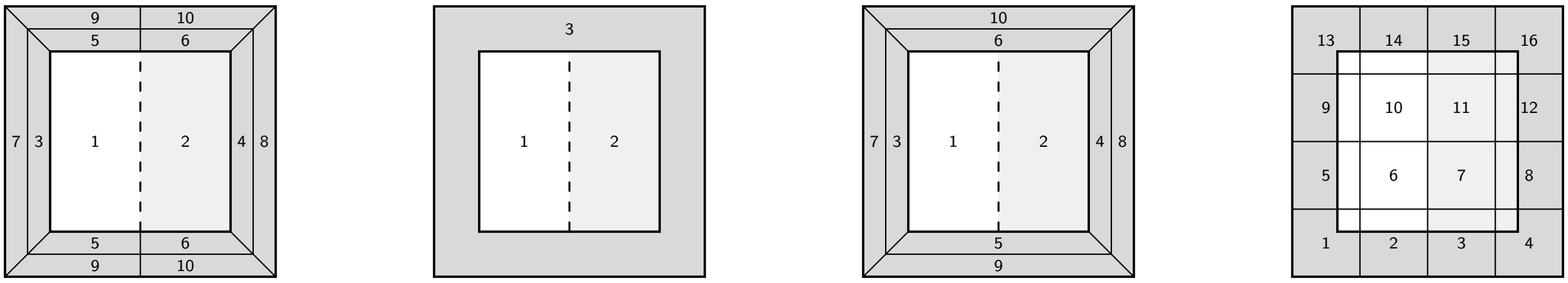
\captionof{figure}{Half-space problem in Section~\ref{sec:helm-correction-rom}, partitions
compared in Table~\ref{tab:reviewer-partition}; (a) is the default partition
of Figure~\ref{fig:halfspace-geometry}.  The PML is dark gray, and pieces with
the same number form one subdomain.}
\label{fig:halfspace-partitions}
\end{table}

\begin{figure}
\centering
\includegraphics[width=0.45\linewidth]{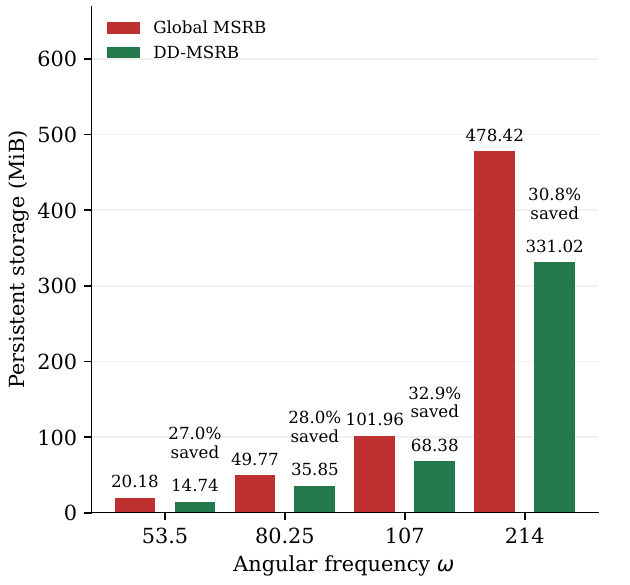}
\includegraphics[width=0.45\linewidth]{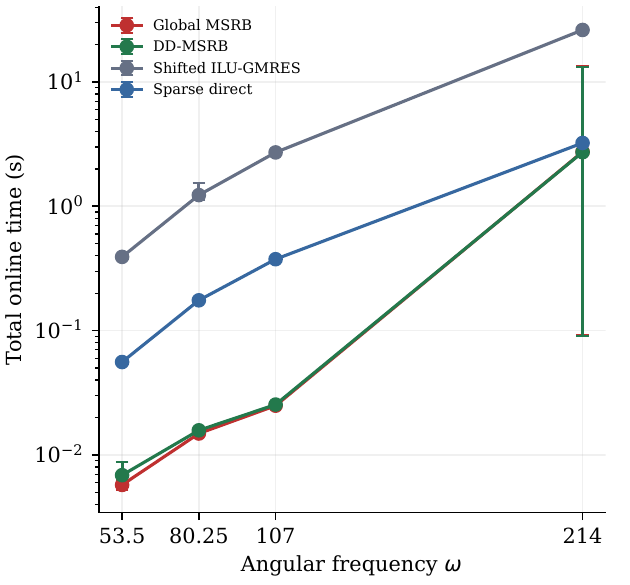}
\caption{Half-space problem in Section~\ref{sec:helm-correction-rom}, frequency study.
Left: ROM payloads and DD-MSRB savings.  Right: mean total online time.}
\label{fig:halfspace-frequency-scaling}
\label{fig:halfspace-frequency-scaling:a}
\end{figure}

\experimentlabel{Effect of frequency}
We increase \(\omega\) through \(53.5\), \(80.25\), \(107\), and \(214\) at
fixed points per wavelength, with approximately constant physical PML thickness
and \(\sigma_{\max}/\omega=80/53.5\).  The extended system grows from
\(35{,}721\) to \(570{,}025\) unknowns. Figure~\ref{fig:halfspace-geometry}
(bottom right) shows a reference solution at \(\omega=214\).  At each frequency, the
models are trained independently on the same 24 values of \(c_{\rm right}\) using the ten-subdomain
partition and \(\eta_j=10^{-5}\), and tested at \(c_{\rm right}=0.8053\),
\(0.8297\), \(0.8761\), \(0.9229\), and \(0.9431\) with tolerance
\(\rho\le10^{-10}\).  Figure~\ref{fig:halfspace-frequency-scaling} shows that
DD maintains a payload reduction of \(27.0\%\)--\(32.9\%\) across this range,
saving \(147.40\) MiB at \(\omega=214\), while its online times remain
comparable to those of global MSRB.  All five queries converge within the ROM
steps at the first three frequencies.  At \(\omega=214\), one query requires
additional SL steps (14 for global MSRB and 13 for DD-MSRB), raising the mean online
times to \(2.7452\) and \(2.7270\) s, respectively.  Thus, the storage benefit
persists at higher frequencies, although occasional SL corrections can increase the mean online cost due to the preconditioner setup.

\subsection{Experiment 2: Elliptic-obstacle scattering}
\label{sec:eccentric}
\experimentlabel{Problem setting}
This experiment tests whether MSRB remains effective for a geometric parameter,
on which the discrete operator depends nonaffinely.  No hyper-reduction is
applied, so the reduced operators are projected online for each test parameter.  We consider a sound-soft ellipse scatterer centered at
\((0.5,0.5)\) in \(\Omega=(0,1)^2\), whose aspect ratio \(\mu=\alpha\) sets the
area-preserving semi-axes \(a=0.12\sqrt{\alpha}\) and \(b=0.12/\sqrt{\alpha}\);
the scattered field satisfies \(u^{\rm sc}=-\exp(\ii\omega x)\) on its
boundary, with \(\omega=18\) and unit background wavespeed.  A second-order
ghost-point discretization with 85 physical points per side and 30 PML cells
per face (\(\sigma_{\max}=80\)) gives \(N_h=21{,}025\).  All geometries share
one background grid: nodes inside the obstacle carry unit matrix rows and
zero right-hand sides and snapshots, the DD index map covers the union of the
active nodes over the training set, and residuals are evaluated on all nodes.
Training uses 96 aspect ratios distributed log-uniformly in \([9/10,10/9]\)
and adjusted by a geometry-admissibility rule; the 22 test geometries differ
from the training geometries.  All methods use the tolerance
\(\rho\le10^{-12}\), and solution errors are measured against sparse direct
solutions.

\begin{figure}
\centering
\includegraphics[width=0.45\linewidth]{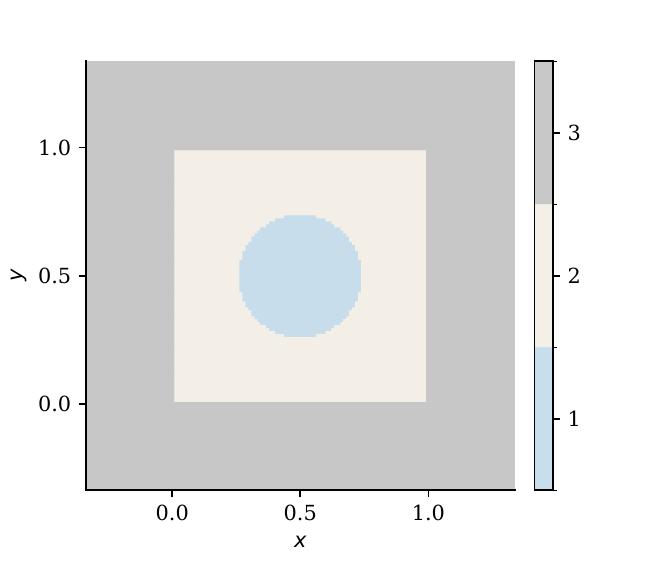}
\includegraphics[width=0.45\linewidth]{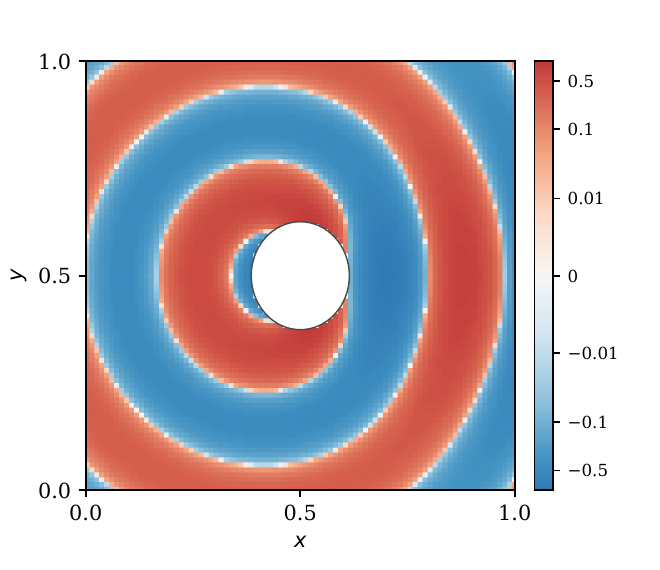}
\caption{Elliptic-obstacle scattering problem in Section~\ref{sec:eccentric}.  Left:
geometry and DD partition into the near-field disk (subdomain 1), the
remaining physical domain (subdomain 2), and the PML (subdomain 3).  Right:
normalized real part of the reference solution in the physical domain at
\(\alpha=0.9092\).}
\label{fig:ecc-setting}
\end{figure}

\begin{figure}
\centering
\includegraphics[width=0.45\linewidth]{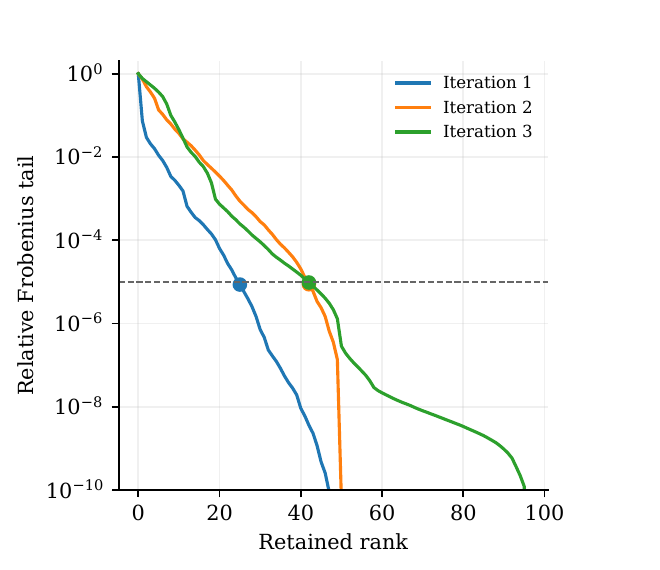}
\includegraphics[width=0.45\linewidth]{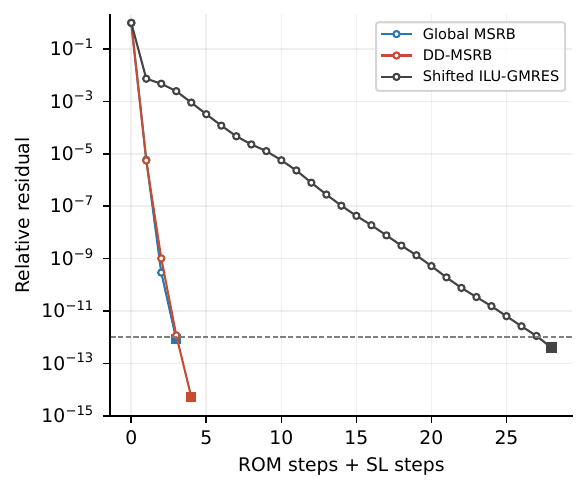}
\caption{Elliptic-obstacle scattering problem in Section~\ref{sec:eccentric}.  Left:
relative Frobenius tails of the iteration-wise global snapshots.  Right:
relative residual histories at \(\alpha=0.9092\).}
\label{fig:ecc-residual-spectra}
\end{figure}

The DD partition separates a near-field disk of radius \(0.24\) around the
obstacle from the remaining physical region and from the PML.  The disk contains the obstacle
at every aspect ratio and confines the part of the field that moves with the
boundary, so one partition serves all geometries, and the far field and the PML
can be represented with lower local ranks.  Figure~\ref{fig:ecc-setting} shows the domain partition (left) and the normalized real part of the solution at the representative query \(\alpha=0.9092\) (right).  Both global MSRB and DD-MSRB are constructed with \(\eta_j=10^{-5}\) for the first three iterations, guided by the relative Frobenius tails of the iteration-wise global
snapshots (Figure~\ref{fig:ecc-residual-spectra}). The retained global ranks are
\((25,42,42)\).

\begin{figure}
\centering
\includegraphics[width=0.45\linewidth]{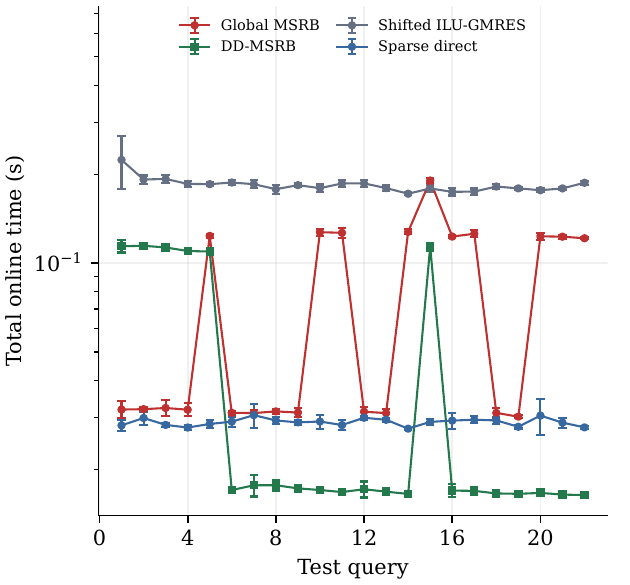}
\hfill
\includegraphics[width=0.45\linewidth]{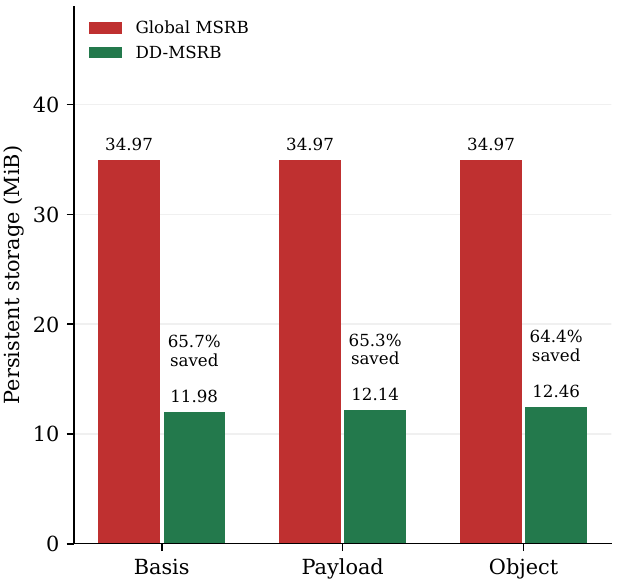}
\caption{Elliptic-obstacle scattering problem in Section~\ref{sec:eccentric},
computational costs.  Left: online times over the 22 test geometries (mean \(\pm\) standard deviation over five repetitions).  Right: total basis storage, deployed
payload, and MATLAB object size.}
\label{fig:ecc-cost}
\end{figure}

\experimentlabel{Convergence and online cost}
MSRB remains effective for this nonaffine geometric parameter.  At
\(\alpha=0.9092\), global MSRB reaches the tolerance within the three ROM steps
and DD-MSRB after one additional SL step, whereas the SL baseline needs 28
iterations (Figure~\ref{fig:ecc-residual-spectra}, right).  Across the 22 test
geometries, global MSRB and DD-MSRB converge within the ROM steps for 12 and 16
geometries and require an average of \(1.14\) and \(0.32\) SL steps, compared
with \(27.68\) iterations for the baseline (Table~\ref{tab:ecc-speedup}).
Their mean online times, \(0.0768\) and \(0.0430\) s, give speedups of
\(2.40\) and \(4.28\) over the baseline's \(0.1840\) s
(Figure~\ref{fig:ecc-cost}, left).  DD-MSRB takes \(44.0\%\) less time than
global MSRB because its ROM steps are cheaper and it needs the ILU setup less
often.  The sparse direct solver remains fastest, at \(0.0290\) s, because this
problem has only \(21{,}025\) unknowns; hyper-reduction
techniques~\cite{BARRAULT2004667,Chaturantabut2010DEIM} or operator
inference~\cite{peherstorfer2016data} could further reduce the online cost of
the ROM steps and are left for future work.  All solutions meet the tolerance,
with maximum relative errors of \(1.19\times10^{-8}\) for global MSRB and
\(9.72\times10^{-9}\) for DD-MSRB.

\begin{table}[!htbp]\centering\normalsize\setlength{\tabcolsep}{3pt}
\setlength{\belowcaptionskip}{6pt}
\caption{Elliptic-obstacle scattering problem in Section~\ref{sec:eccentric}, online
cost over 22 fixed test geometries.  Means include conditional shifted-matrix
construction and ILU factorization (Setup) and the SL steps or the direct
solution (Solve); Total also includes the ROM steps and the final true-residual
check.  Speedup is
relative to Shifted ILU-GMRES.}
\label{tab:ecc-speedup}
\begin{tabular}{lrrrrrr}\toprule
Method & ROM (s) & Setup (s) & Solve (s) & Total (s) & It. & Speedup \\
\midrule
Global MSRB & 0.0306 & 0.0412 & 0.0045 & 0.0768 & 1.14 & $ 2.40 \times $ \\
DD-MSRB & 0.0167 & 0.0249 & 0.0010 & 0.0430 & 0.32 & $ 4.28 \times $ \\
Shifted ILU-GMRES & 0.0000 & 0.0909 & 0.0928 & 0.1840 & 27.68 & $ 1.00 \times $ \\
Sparse direct & 0.0000 & 0.0000 & 0.0287 & 0.0290 & -- & $ 6.35 \times $ \\
\bottomrule\end{tabular}
\end{table}

\experimentlabel{Storage saving}
DD-MSRB saves about two thirds of the ROM storage.  It reduces the basis storage
by \(65.7\%\) (between \(62.0\%\) and \(70.0\%\) per iteration), the payload by
\(65.3\%\), and the MATLAB object size by \(64.4\%\)
(Figure~\ref{fig:ecc-cost}, right).  Its reduced dimensions \((41,71,75)\)
exceed the global ranks \((25,42,42)\), but each local mode is stored only
inside its subdomain.  Because the reduced matrices are assembled online, the
payload holds only the bases and one shared index map.  Offline training takes
about the same time, \(7.57\) s against \(7.60\) s for global MSRB.

\FloatBarrier

\subsection{Experiment 3: Three-dimensional Gaussian waveguide}
\label{sec:3d-waveguide}
To evaluate the proposed framework at larger scale, we next consider a 
three-dimensional Gaussian waveguide at frequencies \(\omega=24\), \(36\), \(48\), and \(60\). At \(\omega=48\), the computational domain spans about 23 wavelengths. In this regime, DD-MSRB requires less
storage and less online time than global MSRB.  We first analyze the
\(\omega=48\) benchmark in detail, then study frequency scaling, and finally
analyze this dual advantage.

\begin{figure}
\centering
\includegraphics[width=0.45\linewidth]{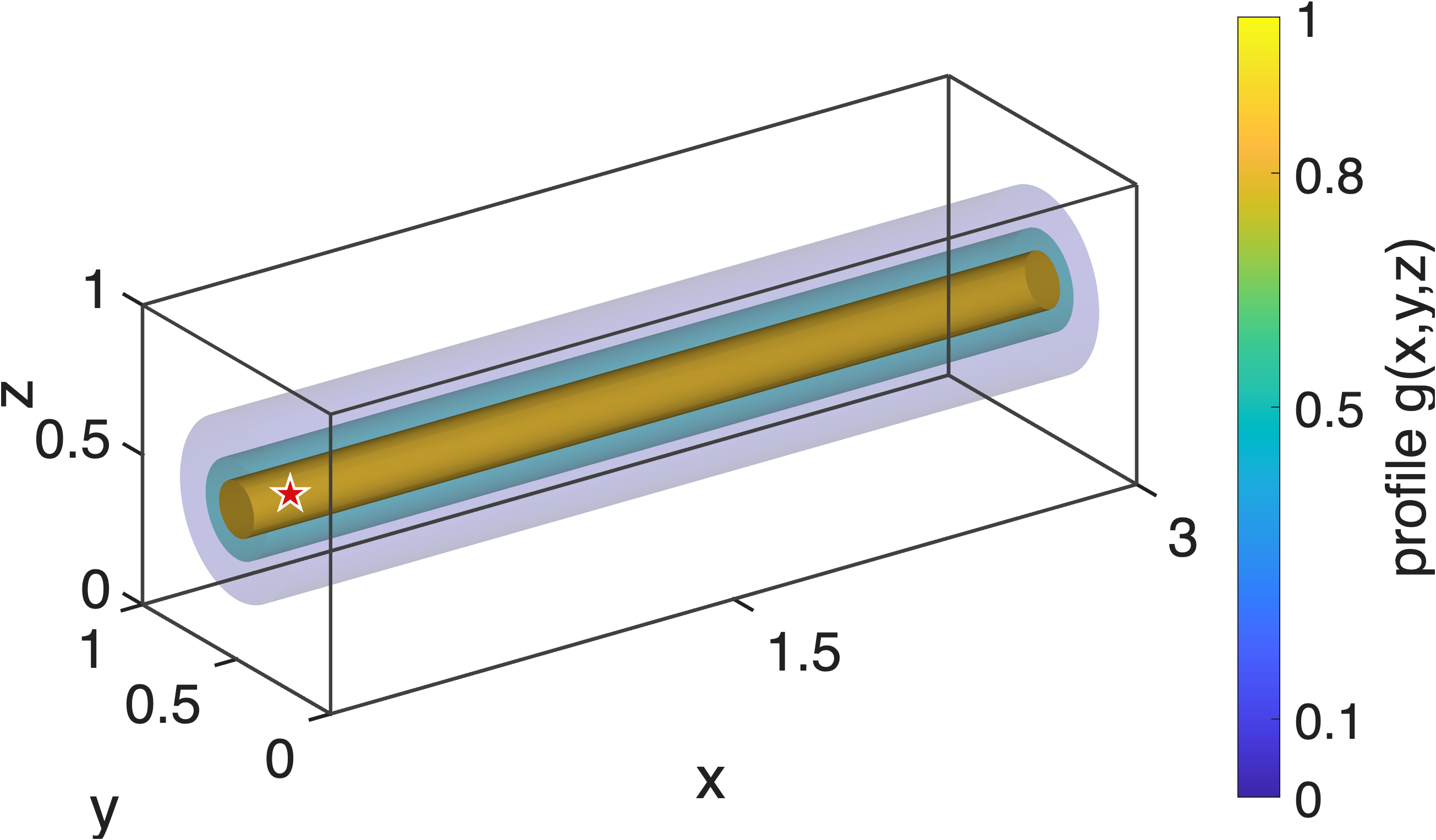}
\hfill
\includegraphics[width=0.45\linewidth]{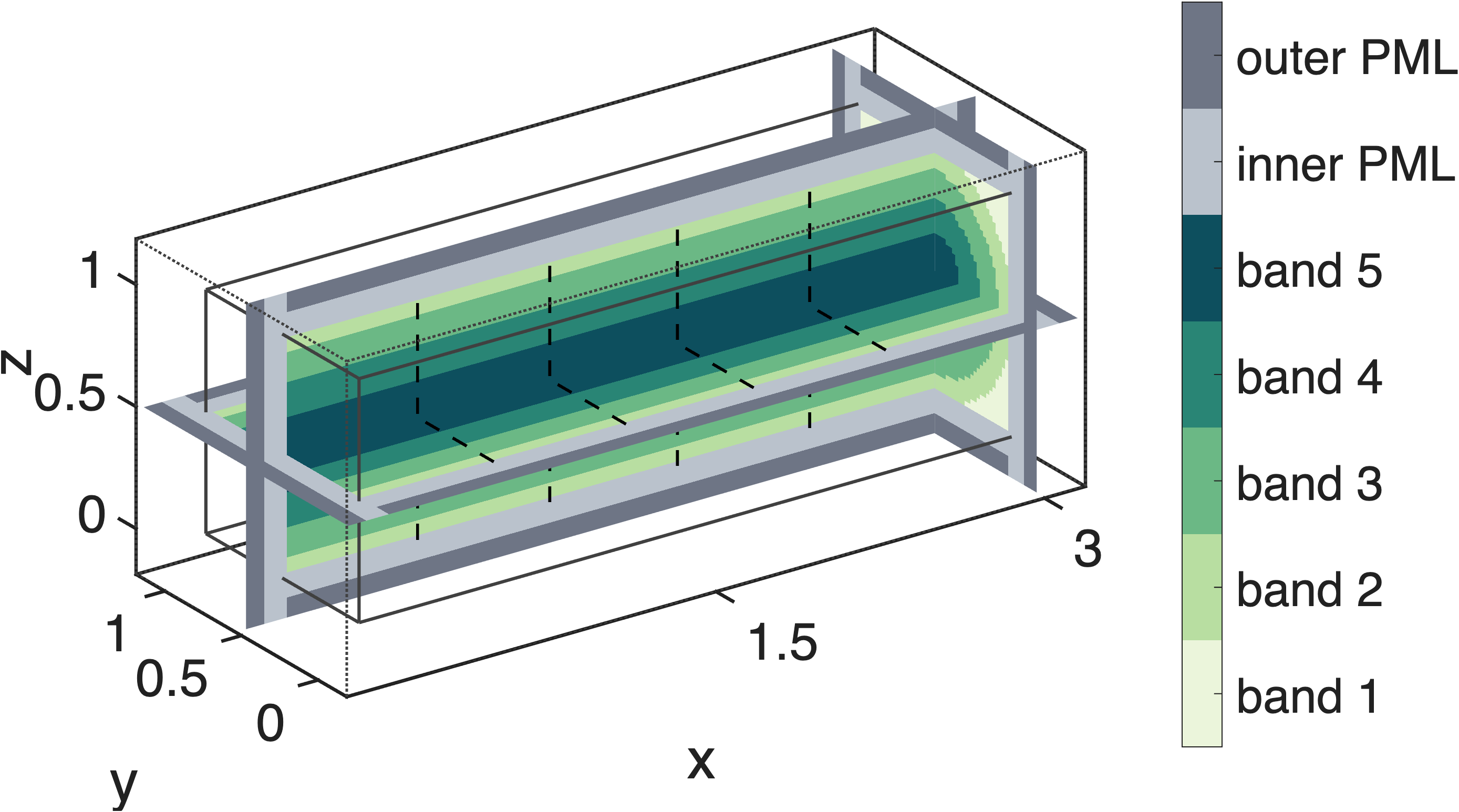}
\caption{Gaussian waveguide problem in Section~\ref{sec:3d-waveguide}.  Left:
guiding-factor isosurfaces and source.  Right: DD partition on longitudinal
and transverse cuts.}
\label{fig:3d-waveguide-geometry}
\end{figure}

\experimentlabel{Problem setting}
The domain is \(\Omega=(0,3)\times(0,1)^2\), with wavespeed profile
\begin{equation}
 c(\vect x;\mu)^{-2}
 =1+\mu\exp\!\left[-\frac12\left(
 \frac{(y-0.50)^2}{0.14^2}+
 \frac{(z-0.50)^2}{0.14^2}\right)\right],
 \label{eq:3d-waveguide-profile}
\end{equation}
and a Gaussian source of standard deviation \(1.5h\) centered at
\((0.20,0.50,0.50)\), where \(h=1/(2\omega)\) is the grid spacing, giving about
ten points per wavelength in the guiding core.  The PML has a fixed width of
\(0.1875\) on every face, rounded to a whole number of cells (9, 14, 18, and 23
cells at \(\omega=24\), 36, 48, and 60), and \(\sigma_{\max}=80\).  Let \(g(y,z)\) denote the
Gaussian factor in~\eqref{eq:3d-waveguide-profile}, so that
\(c^{-2}=1+\mu g\).  The DD partition follows the guiding
profile in the physical domain and the faces of the PML
(Figure~\ref{fig:3d-waveguide-geometry}, right).  A physical node belongs to band
\(\ell\in\{1,\ldots,5\}\) if \(t_{\ell-1}\le g(y,z)<t_\ell\), with
\((t_0,\ldots,t_5)=(0,10^{-3},10^{-2},10^{-1},0.5,\infty)\), and each band is
split into five axial sections of equal length in \(x\), which gives 25
subdomains.  Each PML node is assigned to the first direction, in the order
\(x\), \(y\), \(z\), along which it lies in the PML, and then to the lower or
upper face and to the inner or outer half of the layer, which gives
\(3\times2\times2=12\) subdomains, 37 in total.  The partition thus separates the
guiding core, the exterior, and the absorbing layers.
Figure~\ref{fig:3d-waveguide-geometry} (left) shows the isosurfaces
\(g=0.1,0.5,0.8\), the source as a star, and the inner and outer PML in gray.

At each frequency, the ROMs are built independently with eight uniformly
spaced parameters in \([0.2,0.6]\) and iteration-wise budgets
\((\eta_1,\eta_2,\eta_3)=(0.003,0.008,0.02)\). The budgets grow by a factor of about \(2.5\) per iteration, following the
relaxation principle of inexact Krylov methods: later preconditioning steps acting
on smaller residuals  can be applied less
accurately~\cite{SimonciniSzyld2003}. The SL preconditioner, solver settings, and hardware are described at the
beginning of Section~\ref{sec:numerics}; the SL baseline is labeled Shifted
MG-GMRES in the tables and figures.  At each frequency, we test the ROMs at the same three parameter values:
$\mu=0.237$, $0.317$, and $0.483$. To assess performance beyond these fixed test points,
the detailed study also includes seven additional test parameters
randomly sampled from the intervals between consecutive training
parameters.

\subsubsection{Detailed study at \(\omega=48\)}
\label{sec:3d-representative}
At \(\omega=48\), the \(287\times95\times95\) physical grid padded with an
18-cell PML per face yields \(N_h=5{,}543{,}003\) unknowns.  We test ten
parameters as described above.  All methods for a
given parameter run on the same node (Hygon 7285 or 7285H), and overall speedups are
geometric means of the ten within-node per-parameter ratios.

\begin{figure}
\centering
\begin{minipage}[c]{0.56\linewidth}
\includegraphics[width=\linewidth,trim=68 14 22 31,clip]{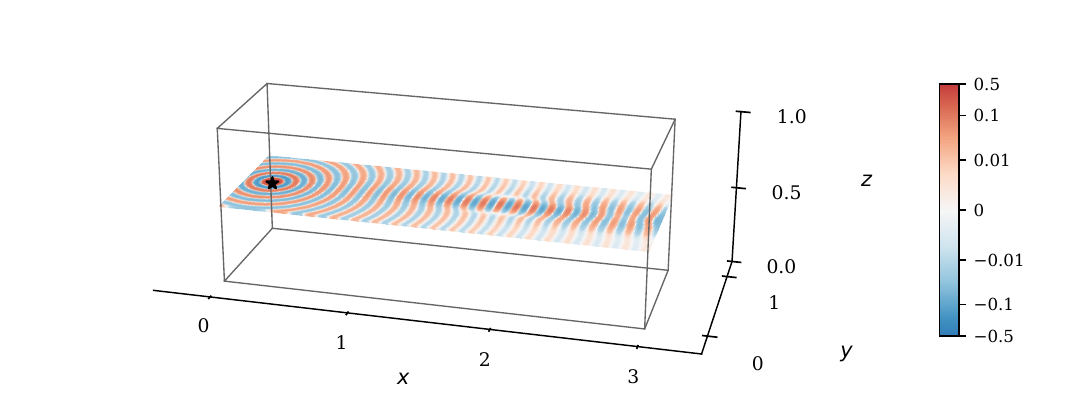}
\end{minipage}
\hfill
\begin{minipage}[c]{0.40\linewidth}
\includegraphics[width=\linewidth]{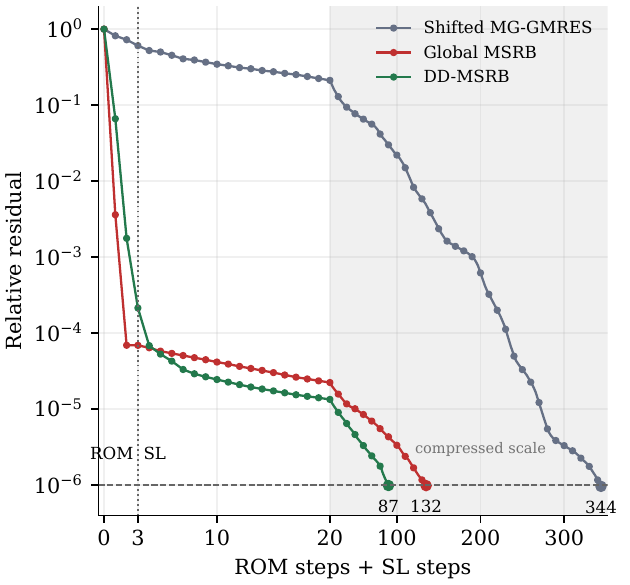}
\end{minipage}
\caption{Gaussian waveguide problem in Section~\ref{sec:3d-representative} at
\(\omega=48\) and \(\mu=0.317\).  Left: normalized real part of the DD-MSRB
solution on \(z=0.5\) in the physical domain.  Right: relative residual
histories over the ROM and SL steps, with the SL step counts marked; the axis
is compressed beyond step 20.}
\label{fig:3d-waveguide-field}
\label{fig:3d-waveguide-residual}
\end{figure}

\experimentlabel{Convergence and online cost}
Figure~\ref{fig:3d-waveguide-field} shows the DD-MSRB solution at the
representative query \(\mu=0.317\) (left) and the residual histories over the
ROM and SL steps (right).  For both MSRB variants, the
three ROM steps from the zero initial guess reduce the relative residual by
nearly four orders of magnitude.  Global MSRB leaves a smaller residual than
DD-MSRB, \(6.93\times10^{-5}\) against \(2.15\times10^{-4}\), but the DD
residual decreases faster in the first SL steps and falls below the global one
(Figure~\ref{fig:3d-waveguide-residual}, right).  DD-MSRB therefore reaches the
tolerance \(10^{-6}\) after 87 SL steps, compared with 132 for global MSRB and
344 iterations for the SL baseline.  The ROM steps take \(3.90\) s for DD-MSRB
and \(3.49\) s for global MSRB, whereas the SL stages take \(240.97\) and
\(366.96\) s (means of three repetitions); Section~\ref{sec:3d-dual-advantage}
analyzes why DD-MSRB needs fewer SL steps.

\begin{figure}
\centering
\includegraphics[width=0.45\linewidth]{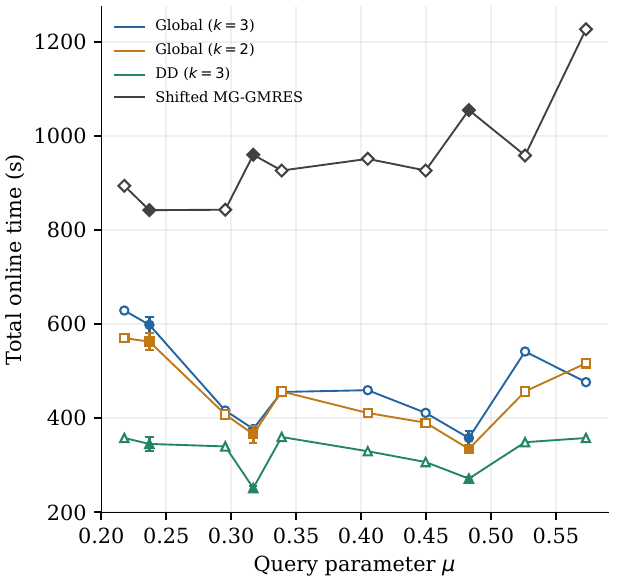}
\hfill
\includegraphics[width=0.45\linewidth]{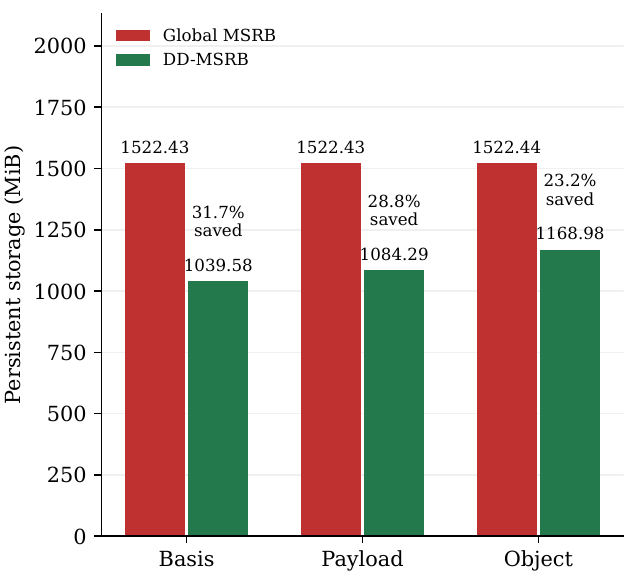}
\caption{Gaussian waveguide problem in Section~\ref{sec:3d-representative} at
\(\omega=48\), computational costs.  Left: total online time at each of the
ten test parameters; filled markers denote the three common queries, with
error bars showing the standard deviation over three repetitions.  Right:
total basis storage, deployed payload, and MATLAB object size.}
\label{fig:3d-validation}
\end{figure}

Over the ten tests (Table~\ref{tab:3d-validation} and
Figure~\ref{fig:3d-validation}, left), DD-MSRB and global MSRB require an average of
\(118.8\) and \(173.5\) SL steps, compared with \(356.2\) iterations for the SL
baseline.  The mean online times are \(326.63\) s for DD-MSRB, \(472.08\) s for
global MSRB, and \(958.19\) s for the baseline.  DD-MSRB is faster than global
MSRB at every test, with geometric-mean speedups of \(1.43\times\) over global
MSRB and \(2.94\times\) over the baseline.

To assess whether the speedup of DD-MSRB also holds at comparable storage,
we test global MSRB with \(k=2\). Omitting the third global reduced
space also reduces the average SL step count from \(173.5\) to
\(168.7\). This variant uses \(1099.54\) MiB, close to the
\(1084.29\) MiB of DD-MSRB, and takes \(447.04\) s on average.
DD-MSRB thus remains \(1.36\times\) faster at comparable storage
(Table~\ref{tab:3d-validation}).

\begin{table}[!htbp]\centering\normalsize
\setlength{\belowcaptionskip}{6pt}
\caption{Gaussian waveguide problem in Section~\ref{sec:3d-representative} at
\(\omega=48\), online cost over the ten test parameters by query group.  Each
group weights its queries equally, and the times exclude matrix assembly and
model loading.}
\label{tab:3d-validation}\par\medskip
\setlength{\tabcolsep}{3pt}
\begin{tabular*}{\linewidth}{@{\extracolsep{\fill}}lrrrrrr@{}}\toprule
& \multicolumn{2}{c}{Three common queries} & \multicolumn{2}{c}{Seven added queries} & \multicolumn{2}{c}{All ten queries} \\
\cmidrule(lr){2-3}\cmidrule(lr){4-5}\cmidrule(l){6-7}
Method & Online (s) & SL steps & Online (s) & SL steps & Online (s) & SL steps \\
\midrule
Global (\(k=3\)) & 444.24 & 155.3 & 484.02 & 181.3 & 472.08 & 173.5 \\
Global (\(k=2\)) & 420.85 & 149.3 & 458.26 & 177.0 & 447.04 & 168.7 \\
DD (\(k=3\)) & 289.11 & 100.7 & 342.71 & 126.6 & 326.63 & 118.8 \\
Shifted MG-GMRES & 952.18 & 349.7 & 960.76 & 359.0 & 958.19 & 356.2 \\
\bottomrule\end{tabular*}
\end{table}

\experimentlabel{Storage saving}
DD-MSRB reduces the basis storage by \(31.7\%\), the payload by \(28.8\%\), and
the MATLAB object size by \(23.2\%\) (Figure~\ref{fig:3d-validation}, right).
The saving concentrates in the first two iterations, \(59.9\%\) and \(28.0\%\).
The 37 local ranks sum to 112, 167, and 197 in the three iterations, against
global ranks of 7, 6, and 5; each local vector is supported on a single
subdomain, so the larger local ranks still take less storage.  In the third
iteration, however, DD-MSRB stores \(3.2\%\) more than global MSRB.  The aggregate time of all
training solves is \(24{,}340\) s for DD-MSRB and \(27{,}033\) s for global MSRB
(Table~\ref{tab:3d-frequency}).

\FloatBarrier
\experimentlabel{Break-even point}
DD-MSRB reduces the mean online time per parameter from \(958.19\) s to
\(326.63\) s, saving \(631.56\) s per parameter relative to the SL
baseline. At this average saving, its aggregate training-solve
cost of \(24{,}340\) s is approximately recovered after
\(  Q_{\mathrm{BE}}
  = \left\lceil \frac{24{,}340}{958.19-326.63} \right\rceil
  = 39\)
tests. This estimate accounts for training solves and online
solution times, excluding other offline costs, matrix assembly,
and model loading.


\subsubsection{Frequency scaling}
\label{sec:3d-frequency}
We repeat the comparison at \(\omega=24\), \(36\), \(48\), and \(60\) with the
discretization and PML of the problem setting, so the system size grows from
\(680{,}225\) at \(\omega=24\) to \(11{,}026{,}125\) at \(\omega=60\).  The runs at \(\omega=24\), \(36\), and \(48\) share one node, while runs at
\(\omega=60\) use another node of the same processor type; at \(\omega=48\), the
online times are the repeated timings of Section~\ref{sec:3d-representative}.  


\begin{table}[!htbp]\centering\small
\setlength{\belowcaptionskip}{6pt}
\caption{Gaussian waveguide problem in Section~\ref{sec:3d-frequency}, frequency study:
relative residual after the ROM steps \(\rho_{\rm ROM}\), SL steps, and total
online time, each averaged over the three common queries, and aggregate time
of all training solves.  Baseline: Shifted MG-GMRES.}
\label{tab:3d-frequency}
\setlength{\tabcolsep}{2pt}
\begin{tabular*}{\linewidth}{@{\extracolsep{\fill}}rrrrrrrrrrr@{}}\toprule
& \multicolumn{2}{c}{\(\rho_{\rm ROM}\)} & \multicolumn{3}{c}{SL steps} & \multicolumn{3}{c}{Online (s)} & \multicolumn{2}{c}{Offline (s)} \\
\cmidrule(lr){2-3}\cmidrule(lr){4-6}\cmidrule(lr){7-9}\cmidrule(l){10-11}
$\omega$ & Global & DD & Global & DD & Baseline & Global & DD & Baseline & Global & DD \\
\midrule
24 & \(1.1\times10^{-5}\) & \(9.2\times10^{-5}\) & 41.3 & 72.0 & 365.3 & 20.13 & 33.38 & 159.86 & 4928 & 4405 \\
36 & \(1.8\times10^{-4}\) & \(2.3\times10^{-4}\) & 110.0 & 84.0 & 324.0 & 128.64 & 97.35 & 375.45 & 14108 & 12537 \\
48 & \(4.4\times10^{-4}\) & \(5.4\times10^{-4}\) & 155.3 & 100.7 & 349.7 & 444.24 & 289.11 & 952.18 & 27033 & 24340 \\
60 & \(4.8\times10^{-3}\) & \(8.0\times10^{-4}\) & 525.0 & 366.3 & 994.3 & 3016.59 & 2106.98 & 5718.79 & 180126 & 161738 \\
\bottomrule\end{tabular*}
\end{table}

\begin{figure}
\centering
\includegraphics[width=0.45\linewidth]{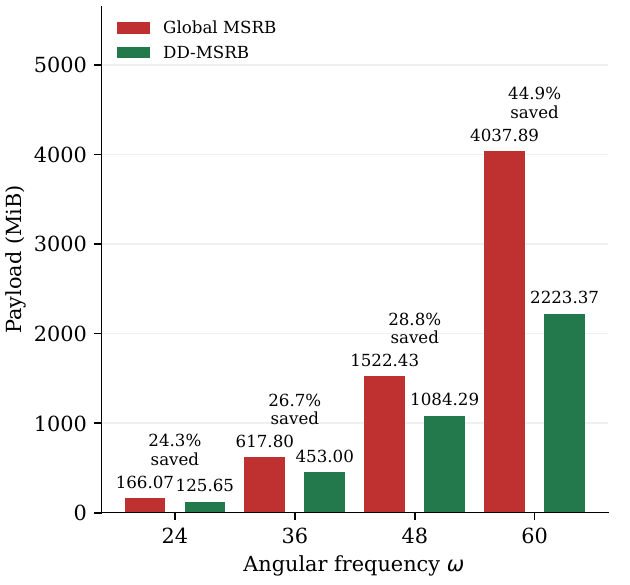}
\includegraphics[width=0.45\linewidth]{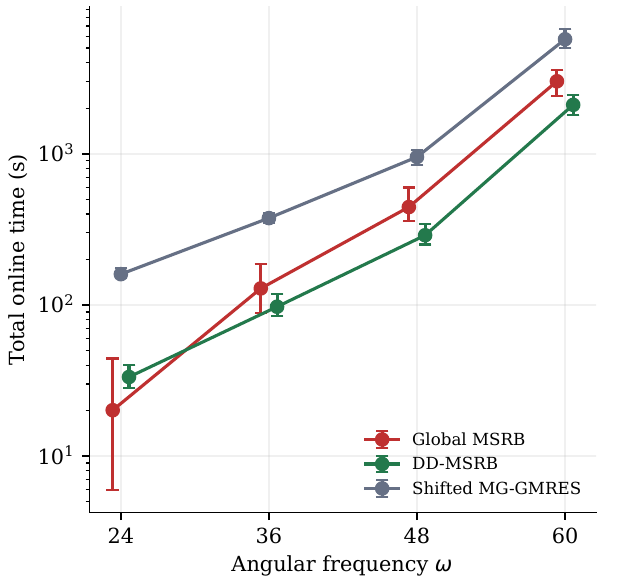}
\caption{Gaussian waveguide problem in Section~\ref{sec:3d-frequency}, frequency study.
Left: ROM payloads and DD-MSRB savings.  Right: mean total online time; error
bars span the three common queries.}
\label{fig:3d-frequency-study}
\end{figure}

\experimentlabel{Storage and online time savings}
The payload saving of DD-MSRB grows with frequency, from \(24.3\%\) at
\(\omega=24\) to \(44.9\%\) at \(\omega=60\)
(Figure~\ref{fig:3d-frequency-study}, left), and the saving in MATLAB object
size from \(18.0\%\) to \(40.8\%\).  At \(\omega=60\), the dimensions of the global MSRB spaces are all \(8\), equal to the size of the training set.  The online time increasingly favors DD-MSRB as the
frequency grows (Table~\ref{tab:3d-frequency} and
Figure~\ref{fig:3d-frequency-study}, right).  At \(\omega=24\), global MSRB is faster on average.  At
\(\omega=36\), DD-MSRB takes \(24.3\%\) less time on average, although it is
slightly slower at one test parameter, where both methods take 77 SL steps.  At
\(\omega=48\) and \(60\), DD-MSRB is faster at every test parameter, with mean savings
of \(34.9\%\) and \(30.2\%\).
With eight training parameters and the offline times of
Table~\ref{tab:3d-frequency}, DD-MSRB recovers its offline cost against the SL baseline after \(35\), \(46\), \(37\), and
\(45\) parameters at \(\omega=24\), \(36\), \(48\), and \(60\), compared with
\(36\), \(58\), \(54\), and \(67\) for global MSRB.

\subsubsection{Dual storage--time advantage of DD-MSRB}
\label{sec:3d-dual-advantage}
For \(\omega\ge36\), DD-MSRB requires less storage and less mean online time
than global MSRB (Section~\ref{sec:3d-frequency}).  The time saving comes from
the SL stage: the three ROM steps take at most \(2.1\%\) of the online time, so
the online time follows the number of SL steps (Table~\ref{tab:3d-frequency}).
The hybrid methods share the full-order operator, SL preconditioner, and
stopping tolerance, and enter the SL stage with the residual \(\rho_{\rm ROM}\)
left by their ROM steps.

At \(\omega=48\), DD-MSRB enters the SL stage with a larger residual at the
three shared test parameters, \(1.21\times10^{-3}\), \(2.15\times10^{-4}\), and
\(1.78\times10^{-4}\) against \(1.17\times10^{-3}\), \(6.93\times10^{-5}\), and
\(7.46\times10^{-5}\) for global MSRB, yet requires fewer SL steps. The residual
norm alone therefore does not determine convergence.  At \(\mu=0.317\), the DD
residual falls below the global residual within the first SL steps and stays
below it until convergence (Figure~\ref{fig:3d-waveguide-residual}, right).
The diagnostics in Appendix~\ref{app:dd-residual-mechanism} show that
DD-MSRB leaves a smaller residual norm in the subdomain interiors
than global MSRB, but a larger residual norm near the interfaces.
The SL pre-smoothing correction also provides a more effective
direction for reducing the DD residual than the global residual.
These observations are consistent with the local ROM corrections
reducing error components that SL damps slowly while leaving
interface mismatches that respond more effectively to smoothing.
This offers a possible explanation for the smaller number of SL
steps despite the larger residual after the ROM stage.

At \(\omega=60\), DD-MSRB has both a smaller residual and fewer SL steps,
\(366.3\) on average against \(525.0\) for global MSRB.  Here global MSRB keeps
all eight training directions in each of the three iterations, the most that
eight training parameters allow, whereas the 37 local ranks of DD-MSRB sum to
121, 174, and 195.  From the same training set, the local bases thus retain
more directions, and the DD residual after the ROM steps is
\(1.81\times10^{-3}\), \(2.79\times10^{-4}\), and \(3.12\times10^{-4}\) at the
three test parameters, against \(1.30\times10^{-2}\), \(8.20\times10^{-4}\),
and \(7.53\times10^{-4}\) for global MSRB.
\FloatBarrier

\section{Conclusion}
\label{sec:conclusion}

Repeated parametric Helmholtz--PML solves offer opportunities
to reuse information that standalone shifted-Laplacian
solves do not exploit. Trajectory-aware MSRB preconditioners
capture this shared structure, but storing a global basis
for each iteration can require substantial memory.
We have developed DD-MSRB to reduce this storage cost
by adapting local basis ranks to the different compressibility
of the physical domain and the PML. The local spaces are
constructed iteration by iteration under a global
error budget through a greedy rank-selection strategy. The resulting ROM preconditioner is further
combined with shifted-Laplacian preconditioning to damp the error components not well represented by the ROMs. 

Across experiments with material and geometric parameters in two and three
dimensions, MSRB shows acceleration over the shifted-Laplacian baseline in
every setting.  DD-MSRB further reduces the ROM storage in almost all settings
tested.  In the three-dimensional waveguide, DD-MSRB also saves online time over global
MSRB as the problem size grows.  We interpret this as a sign that
the domain decomposition is more than a compressed representation.  Its
local structures may more effectively remove the components that the
shifted-Laplacian steps reduce slowly, at the cost of an interface mismatch
that multigrid smoothing damps quickly.

Two limitations of the present work suggest directions for future research.
First, the partitions are chosen from the physical structure of each problem.
An error-guided automatic partitioning strategy similar to~\cite{carlberg2015adaptive} could improve the robustness of the proposed method.  Second, the training parameters are preselected. Adaptive or greedy sampling may be leveraged to further reduce the offline time.
Finally, we plan to apply DD-MSRB to uncertainty quantification and PDE-constrained
optimization, where many forward solves can amortize the offline cost.

\appendix

\section{Residual composition and multigrid response}
\label{app:dd-residual-mechanism}
At \(\omega=48\), the ROM-preconditioned steps of DD-MSRB leave a larger residual than those of
global MSRB, yet fewer SL steps follow (Sections~\ref{sec:3d-representative}
and~\ref{sec:3d-dual-advantage}).  We interpret this as a trade-off.  The subdomain-local
coefficients of DD more effectively remove the error components that the SL preconditioner damps
slowly, but may lead to a larger interface mismatch. Such interface mismatches are high-frequency and can be damped effectively by the multigrid SL preconditioner.  At $\mu=0.317$, the DD residual falls below the global
residual within the first few SL steps and remains
lower until convergence (Figure~\ref{fig:3d-waveguide-residual}, right).
We examine this behavior using the three-iteration ROMs
from the four-frequency study in Section~\ref{sec:3d-frequency}.

The local corrections can be discontinuous across subdomain
interfaces, introducing residual contributions at nodes whose
finite-difference stencils cross an interface. At $\mu=0.317$
and $\omega=36,48,60$, these nodes account for only
$13.2\%$--$21.6\%$ of the physical-domain nodes but contain
$88.5\%$--$90.2\%$ of the DD residual energy, measured by
the squared Euclidean norm. The corresponding fraction
for global MSRB is $2.3\%$--$7.8\%$.

\begin{figure}[!htbp]
\centering
\includegraphics[width=0.49\linewidth]{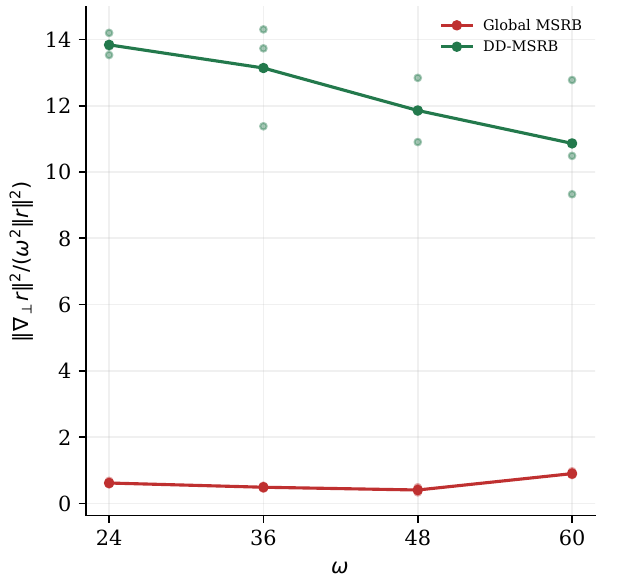}
\hfill
\includegraphics[width=0.49\linewidth]{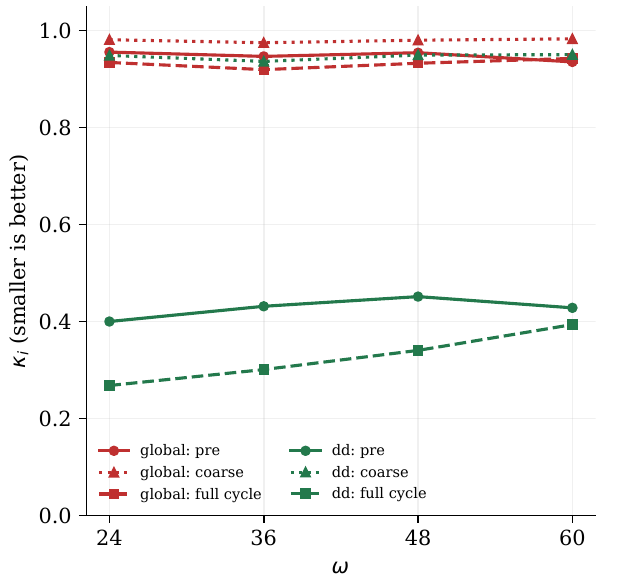}
\caption{Gaussian waveguide problem in Section~\ref{sec:3d-waveguide}, residual and MG
diagnostics.  Left: transverse roughness.  Right: correction-direction
responses \eqref{eq:mg-direction-response}.}
\label{fig:dd-residual-mechanism}
\end{figure}

\begin{figure}[!htbp]
\centering
\includegraphics[width=0.49\linewidth]{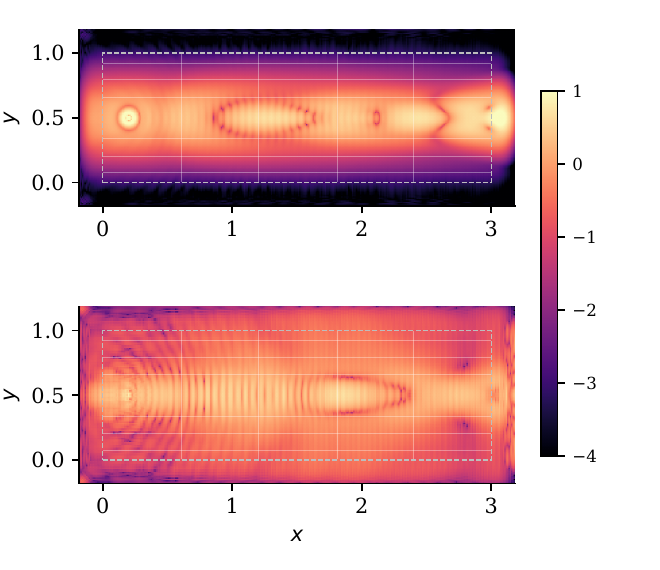}
\hfill
\includegraphics[width=0.49\linewidth]{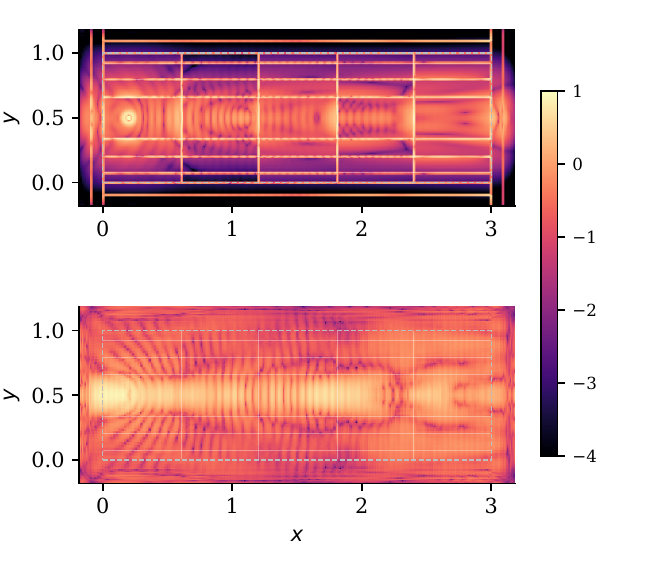}
\caption{Gaussian waveguide problem in Section~\ref{sec:3d-representative} at
\(\omega=48\) and \(\mu=0.317\), central residual slices (\(z=0.5\)).  Left:
global MSRB.  Right: DD-MSRB.  Each side shows the residual at the
ROM-to-SL switch (top) and the residual remaining after 40 SL steps (bottom).}
\label{fig:dd-residual-slices}
\end{figure}

The DD residual also varies more sharply in the transverse
directions. We measure this variation by the transverse
finite-difference energy divided by
$\omega^2\|\vect r_{\phys}\|_2^2$. This quantity is
$10.9$--$13.8$ for DD-MSRB, compared with $0.4$--$0.9$ for global MSRB
(Figure~\ref{fig:dd-residual-mechanism}, left).
The residual slices show the same contrast: the DD residual
is concentrated along subdomain interfaces, while the global
residual is concentrated in the guiding core
(Figure~\ref{fig:dd-residual-slices}, top row).

To examine how multigrid acts on these residuals, we separate
one V-cycle into pre-smoothing, coarse-grid, and post-smoothing
corrections, all represented on the fine grid:
\[
\mathcal B_{\mathrm{SL}}\vect r
=\vect z_{\rm pre}+\vect z_{\rm coarse}+\vect z_{\rm post}.
\]
For each correction, we measure the smallest relative
unshifted residual obtainable by scaling that correction:
\begin{equation}
\kappa_i(\vect r)=\min_{\alpha\in\mathbb C}
\frac{\|\vect r-\alpha\mat A_\mu\vect z_i\|_2}{\|\vect r\|_2},
\qquad i\in\{\mathrm{pre},\mathrm{coarse},\mathrm{post},\mathrm{full}\},
\quad \vect z_{\rm full}=\mathcal B_{\mathrm{SL}}\vect r.
\label{eq:mg-direction-response}
\end{equation}
For DD-MSRB, the mean $\kappa_i$ is $0.40$--$0.45$ for
pre-smoothing, $0.94$--$0.95$ for the coarse-grid correction,
and $0.27$--$0.39$ for the full cycle
(Figure~\ref{fig:dd-residual-mechanism}, right).
Thus, the pre-smoothing correction provides a substantially
more effective direction for reducing the DD residual than
the coarse-grid correction alone.
For global MSRB, all three mean values lie between $0.92$ and $0.98$,
indicating much less effective correction directions.

These results are consistent with DD removing more of the
components that SL reduces slowly, while leaving interface
mismatches that are readily reduced by smoothing.
This interpretation helps explain the rapid initial decrease
of the DD residual despite its larger norm after the ROM stage
(Figure~\ref{fig:3d-waveguide-residual}, right).
A similar residual structure is observed at \(\omega=24\),
but there global MSRB nearly converges within the ROM steps at two of the
three test parameters, and DD-MSRB needs more SL steps on average,
\(72.0\) against \(41.3\) (Table~\ref{tab:3d-frequency}).

\FloatBarrier
\section*{CRediT authorship contribution statement}

\textbf{An Ping:} Conceptualization, Methodology, Software, Investigation, Writing -- original draft.
\textbf{Zhichao Peng:} Conceptualization, Methodology, Supervision, Writing -- review \&
editing.
\textbf{Pengliang Yang:} Methodology, Supervision, Writing -- review \&
editing.
All authors read and approved the final manuscript.

\section*{Declaration of competing interest}

The authors declare that they have no known competing financial interests or
personal relationships that could have appeared to influence the work
reported in this paper.

\section*{Data and code availability}

The code and data will be made publicly available in the GitHub repository
``DD-MSRB preconditioner for Helmholtz''
(\url{https://github.com/frestarcap/DD-MSRB-preconditioner-for-Helmholtz}).
The repository will contain a documented MATLAB implementation of the method with
one driver per experiment, the numerical data behind every table and figure
with the scripts that regenerate them, and the original research code with
which Experiment~3 and Appendix~A were computed, together with reproduction
instructions for each experiment.  All numerical data are generated by the
experiments described in this paper; no external datasets are used.

\section*{Acknowledgements}

An Ping was supported by the Hong Kong PhD Fellowship Scheme (HKPFS). Zhichao Peng was partially supported by the Research Grants Council of Hong Kong (Early Career Scheme 26302724; General Research Fund 16306825 and 16307226).

\section*{Declaration of generative AI and AI-assisted technologies in the manuscript preparation process}

During the preparation of this work, the authors used Claude (Anthropic) and Codex (OpenAI) in order to improve the language and readability of the
manuscript, edit and organize the LaTeX source, check the bibliographic
references, and assist in writing the scripts for the numerical experiments and
figures.  After using these tools, the authors reviewed and edited the content
as needed and take full responsibility for the content of the publication.

\bibliographystyle{elsarticle-num}
\begingroup
\providecommand{\href}[2]{#2}
\let\path\url
\bibliography{references}

\begin{thebibliography}{10}
\expandafter\ifx\csname url\endcsname\relax
  \def\url#1{\texttt{#1}}\fi
\expandafter\ifx\csname urlprefix\endcsname\relax\def\urlprefix{URL }\fi
\expandafter\ifx\csname href\endcsname\relax
  \def\href#1#2{#2} \def\path#1{#1}\fi

\bibitem{BaylissTurkel1980}
A.~Bayliss, E.~Turkel, Radiation boundary conditions for wave-like equations,
  Communications on Pure and Applied Mathematics 33~(6) (1980) 707--725.
\newblock \href {https://doi.org/10.1002/cpa.3160330603}
  {\path{doi:10.1002/cpa.3160330603}}.

\bibitem{Berenger1994}
J.-P. Berenger, A perfectly matched layer for the absorption of electromagnetic
  waves, Journal of Computational Physics 114~(2) (1994) 185--200.
\newblock \href {https://doi.org/10.1006/jcph.1994.1159}
  {\path{doi:10.1006/jcph.1994.1159}}.

\bibitem{BabuskaSauter1997}
I.~Babu{\v{s}}ka, S.~A. Sauter, Is the pollution effect of the {FEM} avoidable
  for the {Helmholtz} equation considering high wave numbers?, SIAM Journal on
  Numerical Analysis 34~(6) (1997) 2392--2423.
\newblock \href {https://doi.org/10.1137/S0036142994269186}
  {\path{doi:10.1137/S0036142994269186}}.

\bibitem{Erlangga2008}
Y.~A. Erlangga, Advances in iterative methods and preconditioners for the
  {Helmholtz} equation, Archives of Computational Methods in Engineering 15~(1)
  (2008) 37--66.
\newblock \href {https://doi.org/10.1007/s11831-007-9013-7}
  {\path{doi:10.1007/s11831-007-9013-7}}.

\bibitem{ErnstGander2012}
O.~G. Ernst, M.~J. Gander, Why it is difficult to solve {Helmholtz} problems
  with classical iterative methods, in: Numerical Analysis of Multiscale
  Problems, Springer, Berlin, Heidelberg, 2012, pp. 325--363.
\newblock \href {https://doi.org/10.1007/978-3-642-22061-6_10}
  {\path{doi:10.1007/978-3-642-22061-6_10}}.

\bibitem{BaylissGoldsteinTurkel1983}
A.~Bayliss, C.~I. Goldstein, E.~Turkel, An iterative method for the {Helmholtz}
  equation, Journal of Computational Physics 49~(3) (1983) 443--457.
\newblock \href {https://doi.org/10.1016/0021-9991(83)90139-0}
  {\path{doi:10.1016/0021-9991(83)90139-0}}.

\bibitem{ErlanggaVuikOosterlee2004}
Y.~A. Erlangga, C.~Vuik, C.~W. Oosterlee, On a class of preconditioners for
  solving the {Helmholtz} equation, Applied Numerical Mathematics 50~(3--4)
  (2004) 409--425.
\newblock \href {https://doi.org/10.1016/j.apnum.2004.01.009}
  {\path{doi:10.1016/j.apnum.2004.01.009}}.

\bibitem{ErlanggaOosterleeVuik2006}
Y.~A. Erlangga, C.~W. Oosterlee, C.~Vuik, A novel multigrid based
  preconditioner for heterogeneous {Helmholtz} problems, SIAM Journal on
  Scientific Computing 27~(4) (2006) 1471--1492.
\newblock \href {https://doi.org/10.1137/040615195}
  {\path{doi:10.1137/040615195}}.

\bibitem{VanGijzenErlanggaVuik2007}
M.~B. van Gijzen, Y.~A. Erlangga, C.~Vuik, Spectral analysis of the discrete
  {Helmholtz} operator preconditioned with a shifted {Laplacian}, SIAM Journal
  on Scientific Computing 29~(5) (2007) 1942--1958.
\newblock \href {https://doi.org/10.1137/060661491}
  {\path{doi:10.1137/060661491}}.

\bibitem{Gander2006}
M.~J. Gander, I.~G. Graham, E.~A. Spence, Applying {GMRES} to the {Helmholtz}
  equation with shifted {Laplacian} preconditioning: what is the largest shift
  for which wavenumber-independent convergence is guaranteed?, Numerische
  Mathematik 131~(3) (2015) 567--614.
\newblock \href {https://doi.org/10.1007/s00211-015-0700-2}
  {\path{doi:10.1007/s00211-015-0700-2}}.

\bibitem{BENAMOU199768}
J.-D. Benamou, B.~Desprès, A domain decomposition method for the {Helmholtz}
  equation and related optimal control problems, Journal of Computational
  Physics 136~(1) (1997) 68--82.
\newblock \href {https://doi.org/10.1006/jcph.1997.5742}
  {\path{doi:10.1006/jcph.1997.5742}}.

\bibitem{farhat2000two}
C.~Farhat, A.~Macedo, M.~Lesoinne, {A two-level domain decomposition method for
  the iterative solution of high frequency exterior Helmholtz problems},
  Numerische Mathematik 85~(2) (2000) 283--308.
\newblock \href {https://doi.org/10.1007/PL00005389}
  {\path{doi:10.1007/PL00005389}}.

\bibitem{GanderMagoulesNataf2002}
M.~J. Gander, F.~Magoul{\`e}s, F.~Nataf, Optimized {Schwarz} methods without
  overlap for the {Helmholtz} equation, SIAM Journal on Scientific Computing
  24~(1) (2002) 38--60.
\newblock \href {https://doi.org/10.1137/S1064827501387012}
  {\path{doi:10.1137/S1064827501387012}}.

\bibitem{GrahamSpenceVainikko2017}
I.~G. Graham, E.~A. Spence, E.~Vainikko, Domain decomposition preconditioning
  for high-frequency {Helmholtz} problems with absorption, Mathematics of
  Computation 86~(307) (2017) 2089--2127.
\newblock \href {https://doi.org/10.1090/mcom/3190}
  {\path{doi:10.1090/mcom/3190}}.

\bibitem{EngquistYing2011}
B.~Engquist, L.~Ying, Sweeping preconditioner for the {Helmholtz} equation:
  Moving perfectly matched layers, Multiscale Modeling \& Simulation 9~(2)
  (2011) 686--710.
\newblock \href {https://doi.org/10.1137/100804644}
  {\path{doi:10.1137/100804644}}.

\bibitem{taus2020sweeps}
M.~Taus, L.~Zepeda-N{\'u}{\~n}ez, R.~J. Hewett, L.~Demanet, {L-Sweeps: A
  scalable, parallel preconditioner for the high-frequency Helmholtz equation},
  Journal of Computational Physics 420 (2020) 109706.
\newblock \href {https://doi.org/10.1016/j.jcp.2020.109706}
  {\path{doi:10.1016/j.jcp.2020.109706}}.

\bibitem{DwarkaVuik2022}
V.~Dwarka, C.~Vuik, Scalable multi-level deflation preconditioning for highly
  indefinite time-harmonic waves, Journal of Computational Physics 469 (2022)
  111327.
\newblock \href {https://doi.org/10.1016/j.jcp.2022.111327}
  {\path{doi:10.1016/j.jcp.2022.111327}}.

\bibitem{bristeau1998controllability}
M.-O. Bristeau, R.~Glowinski, J.~P{\'e}riaux, {Controllability methods for the
  computation of time-periodic solutions; application to scattering}, Journal
  of Computational Physics 147~(2) (1998) 265--292.
\newblock \href {https://doi.org/10.1006/jcph.1998.6044}
  {\path{doi:10.1006/jcph.1998.6044}}.

\bibitem{stolk2021time}
C.~C. Stolk, {A time-domain preconditioner for the Helmholtz equation}, SIAM
  Journal on Scientific Computing 43~(5) (2021) A3469--A3502.
\newblock \href {https://doi.org/10.1137/20M1359997}
  {\path{doi:10.1137/20M1359997}}.

\bibitem{grote2019controllability}
M.~J. Grote, J.~H. Tang, {On controllability methods for the Helmholtz
  equation}, Journal of Computational and Applied Mathematics 358 (2019)
  306--326.
\newblock \href {https://doi.org/10.1016/j.cam.2019.03.016}
  {\path{doi:10.1016/j.cam.2019.03.016}}.

\bibitem{appelo2020waveholtz}
D.~Appel{\"o}, F.~Garcia, O.~Runborg, {WaveHoltz: Iterative solution of the
  Helmholtz equation via the wave equation}, SIAM Journal on Scientific
  Computing 42~(4) (2020) A1950--A1983.
\newblock \href {https://doi.org/10.1137/19M1299062}
  {\path{doi:10.1137/19M1299062}}.

\bibitem{LahayeTangVuik2017}
D.~Lahaye, J.~Tang, K.~Vuik (Eds.), Modern Solvers for {Helmholtz} Problems,
  Geosystems Mathematics, Birkh{\"a}user, Cham, 2017.
\newblock \href {https://doi.org/10.1007/978-3-319-28832-1}
  {\path{doi:10.1007/978-3-319-28832-1}}.

\bibitem{GanderZhang2019}
M.~J. Gander, H.~Zhang, A class of iterative solvers for the {Helmholtz}
  equation: Factorizations, sweeping preconditioners, source transfer, single
  layer potentials, polarized traces, and optimized {Schwarz} methods, SIAM
  Review 61~(1) (2019) 3--76.
\newblock \href {https://doi.org/10.1137/16M109781X}
  {\path{doi:10.1137/16M109781X}}.

\bibitem{GanderZhang2022}
M.~J. Gander, H.~Zhang, {Schwarz} methods by domain truncation, Acta Numerica
  31 (2022) 1--134.
\newblock \href {https://doi.org/10.1017/S0962492922000034}
  {\path{doi:10.1017/S0962492922000034}}.

\bibitem{BennerGugercinWillcox2015}
P.~Benner, S.~Gugercin, K.~Willcox, A survey of projection-based model
  reduction methods for parametric dynamical systems, SIAM Review 57~(4) (2015)
  483--531.
\newblock \href {https://doi.org/10.1137/130932715}
  {\path{doi:10.1137/130932715}}.

\bibitem{SenVeroyHuynhDeparisNguyenPatera2006}
S.~Sen, K.~Veroy, D.~B.~P. Huynh, S.~Deparis, N.~C. Nguyen, A.~T. Patera,
  ``{Natural} norm'' a posteriori error estimators for reduced basis
  approximations, Journal of Computational Physics 217~(1) (2006) 37--62.
\newblock \href {https://doi.org/10.1016/j.jcp.2006.02.012}
  {\path{doi:10.1016/j.jcp.2006.02.012}}.

\bibitem{ChenHesthavenMadayRodriguez2010}
Y.~Chen, J.~S. Hesthaven, Y.~Maday, J.~Rodr{\'i}guez, Certified reduced basis
  methods and output bounds for the harmonic {Maxwell}'s equations, SIAM
  Journal on Scientific Computing 32~(2) (2010) 970--996.
\newblock \href {https://doi.org/10.1137/09075250X}
  {\path{doi:10.1137/09075250X}}.

\bibitem{pomplun2010accelerated}
J.~Pomplun, F.~Schmidt, {Accelerated a posteriori error estimation for the
  reduced basis method with application to 3D electromagnetic scattering
  problems}, SIAM Journal on Scientific Computing 32~(2) (2010) 498--520.
\newblock \href {https://doi.org/10.1137/090760271}
  {\path{doi:10.1137/090760271}}.

\bibitem{hess2013fast}
M.~W. Hess, P.~Benner, {Fast evaluation of time--harmonic Maxwell's equations
  using the reduced basis method}, IEEE Transactions on Microwave Theory and
  Techniques 61~(6) (2013) 2265--2274.
\newblock \href {https://doi.org/10.1109/TMTT.2013.2258167}
  {\path{doi:10.1109/TMTT.2013.2258167}}.

\bibitem{benner2015uncertainty}
P.~Benner, J.~Schneider, {Uncertainty quantification for Maxwell's equations
  using stochastic collocation and model order reduction}, International
  Journal for Uncertainty Quantification 5~(3) (2015) 195--208.
\newblock \href
  {https://doi.org/10.1615/Int.J.UncertaintyQuantification.2015010170}
  {\path{doi:10.1615/Int.J.UncertaintyQuantification.2015010170}}.

\bibitem{dal2017numerical}
N.~Dal~Santo, S.~Deparis, A.~Manzoni, {A numerical investigation of multi space
  reduced basis preconditioners for parametrized elliptic advection-diffusion
  equations}, Communications in Applied and Industrial Mathematics 8~(1) (2017)
  282--297.
\newblock \href {https://doi.org/10.1515/caim-2017-0015}
  {\path{doi:10.1515/caim-2017-0015}}.

\bibitem{DalSanto2018}
N.~Dal~Santo, S.~Deparis, A.~Manzoni, A.~Quarteroni, Multi space reduced basis
  preconditioners for large-scale parametrized {PDEs}, SIAM Journal on
  Scientific Computing 40~(2) (2018) A954--A983.
\newblock \href {https://doi.org/10.1137/16M1089149}
  {\path{doi:10.1137/16M1089149}}.

\bibitem{nguyen2018reduced}
N.-C. Nguyen, Y.~Chen, {Reduced-basis method for the iterative solution of
  parametrized symmetric positive-definite linear systems}, arXiv:1804.06363
  [math.NA] (2018).
\newblock \href {https://doi.org/10.48550/arXiv.1804.06363}
  {\path{doi:10.48550/arXiv.1804.06363}}.

\bibitem{liu2022rom}
J.~Liu, Z.~Wang, {A ROM-accelerated parallel-in-time preconditioner for solving
  all-at-once systems in unsteady convection-diffusion PDEs}, Applied
  Mathematics and Computation 416 (2022) 126750.
\newblock \href {https://doi.org/10.1016/j.amc.2021.126750}
  {\path{doi:10.1016/j.amc.2021.126750}}.

\bibitem{mcclarren2022data}
R.~G. McClarren, T.~S. Haut, {Data-driven acceleration of thermal radiation
  transfer calculations with the dynamic mode decomposition and a sequential
  singular value decomposition}, Journal of Computational Physics 448 (2022)
  110756.
\newblock \href {https://doi.org/10.1016/j.jcp.2021.110756}
  {\path{doi:10.1016/j.jcp.2021.110756}}.

\bibitem{peng2024reduced}
Z.~Peng, {Reduced order model enhanced source iteration with synthetic
  acceleration for parametric radiative transfer equation}, Journal of
  Computational Physics 517 (2024) 113303.
\newblock \href {https://doi.org/10.1016/j.jcp.2024.113303}
  {\path{doi:10.1016/j.jcp.2024.113303}}.

\bibitem{Peng2025}
Z.~Peng, A flexible {GMRES} solver with reduced order model enhanced synthetic
  acceleration preconditioner for parametric radiative transfer equation,
  Journal of Computational Physics 534 (2025) 114004.
\newblock \href {https://doi.org/10.1016/j.jcp.2025.114004}
  {\path{doi:10.1016/j.jcp.2025.114004}}.

\bibitem{TangPeng2025}
N.~Tang, Z.~Peng, Synthetic acceleration preconditioners for parametric
  radiative transfer equations based on trajectory-aware reduced order models,
  arXiv:2509.05001 [math.NA] (2025).
\newblock \href {https://doi.org/10.48550/arXiv.2509.05001}
  {\path{doi:10.48550/arXiv.2509.05001}}.

\bibitem{BuffoniTelibIollo2009}
M.~Buffoni, H.~Telib, A.~Iollo, Iterative methods for model reduction by domain
  decomposition, Computers \& Fluids 38~(6) (2009) 1160--1167.
\newblock \href {https://doi.org/10.1016/j.compfluid.2008.11.008}
  {\path{doi:10.1016/j.compfluid.2008.11.008}}.

\bibitem{BaigesCodinaIdelsohn2013}
J.~Baiges, R.~Codina, S.~Idelsohn, A domain decomposition strategy for reduced
  order models. {Application} to the incompressible {Navier--Stokes} equations,
  Computer Methods in Applied Mechanics and Engineering 267 (2013) 23--42.
\newblock \href {https://doi.org/10.1016/j.cma.2013.08.001}
  {\path{doi:10.1016/j.cma.2013.08.001}}.

\bibitem{HoangChoiCarlberg2021}
C.~Hoang, Y.~Choi, K.~Carlberg, Domain-decomposition least-squares
  {Petrov--Galerkin} ({DD-LSPG}) nonlinear model reduction, Computer Methods in
  Applied Mechanics and Engineering 384 (2021) 113997.
\newblock \href {https://doi.org/10.1016/j.cma.2021.113997}
  {\path{doi:10.1016/j.cma.2021.113997}}.

\bibitem{MadayRonquist2002}
Y.~Maday, E.~M. R{\o}nquist, A reduced-basis element method, Comptes Rendus
  Math{\'e}matique 335~(2) (2002) 195--200.
\newblock \href {https://doi.org/10.1016/S1631-073X(02)02427-5}
  {\path{doi:10.1016/S1631-073X(02)02427-5}}.

\bibitem{HuynhKnezevicPatera2013Complex}
D.~B.~P. Huynh, D.~J. Knezevic, A.~T. Patera, A static condensation reduced
  basis element method: Complex problems, Computer Methods in Applied Mechanics
  and Engineering 259 (2013) 197--216.
\newblock \href {https://doi.org/10.1016/j.cma.2013.02.013}
  {\path{doi:10.1016/j.cma.2013.02.013}}.

\bibitem{EftangPatera2013}
J.~L. Eftang, A.~T. Patera, Port reduction in parametrized component static
  condensation: Approximation and a posteriori error estimation, International
  Journal for Numerical Methods in Engineering 96~(5) (2013) 269--302.
\newblock \href {https://doi.org/10.1002/nme.4543}
  {\path{doi:10.1002/nme.4543}}.

\bibitem{smetana2016optimal}
K.~Smetana, A.~T. Patera, {Optimal local approximation spaces for
  component-based static condensation procedures}, SIAM Journal on Scientific
  Computing 38~(5) (2016) A3318--A3356.
\newblock \href {https://doi.org/10.1137/15M1009603}
  {\path{doi:10.1137/15M1009603}}.

\bibitem{ohlberger2015error}
M.~Ohlberger, F.~Schindler, Error control for the localized reduced basis
  multiscale method with adaptive on-line enrichment, SIAM Journal on
  Scientific Computing 37~(6) (2015) A2865--A2895.
\newblock \href {https://doi.org/10.1137/151003660}
  {\path{doi:10.1137/151003660}}.

\bibitem{carlberg2015adaptive}
K.~Carlberg, {Adaptive h-refinement for reduced-order models}, International
  Journal for Numerical Methods in Engineering 102~(5) (2015) 1192--1210.
\newblock \href {https://doi.org/10.1002/nme.4800}
  {\path{doi:10.1002/nme.4800}}.

\bibitem{BuhrSmetana2018}
A.~Buhr, K.~Smetana, Randomized local model order reduction, SIAM Journal on
  Scientific Computing 40~(4) (2018) A2120--A2151.
\newblock \href {https://doi.org/10.1137/17M1138480}
  {\path{doi:10.1137/17M1138480}}.

\bibitem{BuhrEtAl2021}
A.~Buhr, L.~Iapichino, M.~Ohlberger, S.~Rave, F.~Schindler, K.~Smetana,
  Localized model reduction for parameterized problems, in: P.~Benner,
  S.~Grivet-Talocia, A.~Quarteroni, G.~Rozza, W.~Schilders, L.~M. Silveira
  (Eds.), Model Order Reduction, Volume 2: Snapshot-Based Methods and
  Algorithms, De Gruyter, Berlin, 2021, pp. 245--306.
\newblock \href {https://doi.org/10.1515/9783110671490-006}
  {\path{doi:10.1515/9783110671490-006}}.

\bibitem{CocquetGander2017}
P.-H. Cocquet, M.~J. Gander, How large a shift is needed in the shifted
  {Helmholtz} preconditioner for its effective inversion by multigrid?, SIAM
  Journal on Scientific Computing 39~(2) (2017) A438--A478.
\newblock \href {https://doi.org/10.1137/15M102085X}
  {\path{doi:10.1137/15M102085X}}.

\bibitem{etter2020online}
P.~A. Etter, K.~T. Carlberg, {Online adaptive basis refinement and compression
  for reduced-order models via vector-space sieving}, Computer Methods in
  Applied Mechanics and Engineering 364 (2020) 112931.
\newblock \href {https://doi.org/10.1016/j.cma.2020.112931}
  {\path{doi:10.1016/j.cma.2020.112931}}.

\bibitem{rocha2020adaptive}
I.~Rocha, F.~van~der Meer, L.~J. Sluys, {An adaptive domain-based POD/ECM
  hyper-reduced modeling framework without offline training}, Computer Methods
  in Applied Mechanics and Engineering 358 (2020) 112650.
\newblock \href {https://doi.org/10.1016/j.cma.2019.112650}
  {\path{doi:10.1016/j.cma.2019.112650}}.

\bibitem{ferreira2022adaptivity}
B.~P. Ferreira, F.~M. Andrade~Pires, M.~A. Bessa, {Adaptivity for
  clustering-based reduced-order modeling of localized history-dependent
  phenomena}, Computer Methods in Applied Mechanics and Engineering 393 (2022)
  114726.
\newblock \href {https://doi.org/10.1016/j.cma.2022.114726}
  {\path{doi:10.1016/j.cma.2022.114726}}.

\bibitem{BARRAULT2004667}
M.~Barrault, Y.~Maday, N.~C. Nguyen, A.~T. Patera, An empirical interpolation
  method: application to efficient reduced-basis discretization of partial
  differential equations, Comptes Rendus Math{\'e}matique 339~(9) (2004)
  667--672.
\newblock \href {https://doi.org/10.1016/j.crma.2004.08.006}
  {\path{doi:10.1016/j.crma.2004.08.006}}.

\bibitem{Chaturantabut2010DEIM}
S.~Chaturantabut, D.~C. Sorensen, Nonlinear model reduction via discrete
  empirical interpolation, SIAM Journal on Scientific Computing 32~(5) (2010)
  2737--2764.
\newblock \href {https://doi.org/10.1137/090766498}
  {\path{doi:10.1137/090766498}}.

\bibitem{Saad1993}
Y.~Saad, A flexible inner-outer preconditioned {GMRES} algorithm, SIAM Journal
  on Scientific Computing 14~(2) (1993) 461--469.
\newblock \href {https://doi.org/10.1137/0914028} {\path{doi:10.1137/0914028}}.

\bibitem{peherstorfer2016data}
B.~Peherstorfer, K.~Willcox, {Data-driven operator inference for nonintrusive
  projection-based model reduction}, Computer Methods in Applied Mechanics and
  Engineering 306 (2016) 196--215.
\newblock \href {https://doi.org/10.1016/j.cma.2016.03.025}
  {\path{doi:10.1016/j.cma.2016.03.025}}.

\bibitem{SimonciniSzyld2003}
V.~Simoncini, D.~B. Szyld, Theory of inexact {Krylov} subspace methods and
  applications to scientific computing, SIAM Journal on Scientific Computing
  25~(2) (2003) 454--477.
\newblock \href {https://doi.org/10.1137/S1064827502406415}
  {\path{doi:10.1137/S1064827502406415}}.

\end{thebibliography}
\endgroup
\end{document}